\documentclass[graybox]{svmult}

\usepackage{mathptmx}       
\usepackage{helvet}         
\usepackage{courier}        
\usepackage{type1cm}        
\usepackage{amssymb}
\usepackage{amsmath}

\usepackage{makeidx}         
\usepackage{graphicx}        
\usepackage{multicol}        
\usepackage[bottom]{footmisc}

\spdefaulttheorem{project}{Research Project}{\bf}{}
\newcommand{\resproject}[1]{\begin{svgraybox} \begin{project} #1 \end{project} \end{svgraybox}}

\spdefaulttheorem{cproblem}{Challenge Problem}{\bf}{}
\newcommand{\chproblem}[1]{\begin{cproblem} #1 \end{cproblem}}

\spnewtheorem*{prerequisites}{Suggested prerequisites}{\bf}{\small\em}

\makeindex             

\begin{document}

\title*{Tropical Geometry}
\author{Ralph Morrison}
\institute{Ralph Morrison \at Williams College, Address of Institute, \email{10rem@williams.edu}
}
%
%
\maketitle

\abstract*{Tropical mathematics redefines the rules of arithmetic by replacing addition with taking a maximum, and by replacing multiplication with addition. After briefly discussing linear algebra with these operations, we interpret polynomials with these new operations.  This allows us to define piecewise-linear objects called tropical varieties.  We explore these tropical varieties in two and three dimensions, building up discrete tools for studying them and determining their geometric properties.  We then discuss the relationship between tropical geometry and algebraic geometry, which considers shapes defined by usual polynomial equations.}

\abstract{Tropical mathematics redefines the rules of arithmetic by replacing addition with taking a maximum, and by replacing multiplication with addition. After briefly discussing a tropical version of linear algebra, we study polynomials build with these new operations.  These equations define piecewise-linear geometric objects called tropical varieties.  We explore these tropical varieties in two and three dimensions, building up discrete tools for studying them and determining their geometric properties.  We then discuss the relationship between tropical geometry and algebraic geometry, which considers shapes defined by usual polynomial equations.}

\begin{prerequisites}
We use standard set theory notation (unions, functions, etc.) throughout this chapter. Section 1 draws on terminology and motivation from abstract algebra and linear algebra, but can be understood without them.  Section 2 draws on topics from discrete geometry, although it  is mostly self-contained. Section 3 includes geometry in three dimensions, which uses some notation from a standard course in multivariable calculus.  Section \ref{section:tropicalization} uses ring theory terminology from an abstract algebra course. \end{prerequisites}

\section{Tropical Mathematics}

\begin{figure}[hbt]
\sidecaption
\includegraphics[scale=.65]{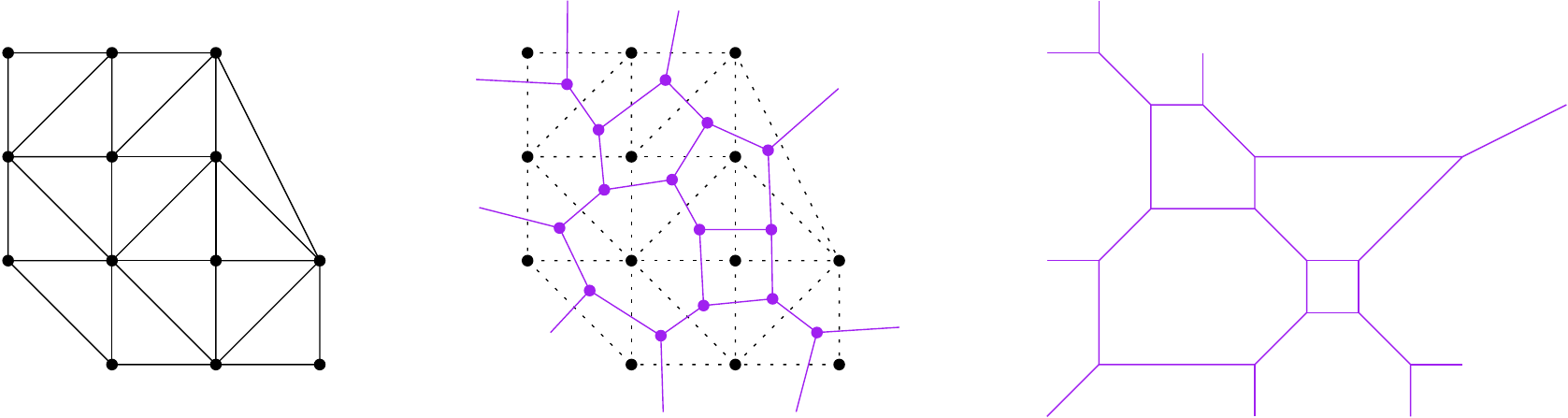}
\caption{Drawing a tropical curve}
\label{fig:drawing}       
\end{figure}

Take a piece of graph paper, or draw your own rectangular grid.  Pick some of the grid points, and join them up to form a polygon.  Be sure it's convex, so that all the angles are less than $180$ degrees.  Now, start connecting grid points to each other with line segments, never letting any two line segments cross.  Keep going until you can't split things up anymore.  You should end up with lots of triangles, like the first picture in Figure \ref{fig:drawing}.

Using a different color, say purple, put a dot in every triangle.  Connect two dots with a line segment if their triangles share a side.  If a triangle has a boundary edge, just draw a little edge coming out of the dot.  Your picture will now look like the middle of Figure \ref{fig:drawing}.  Now, try to draw your purple shape again, but with the following rule:  each line segment you draw should be perpendicular to the shared side of the triangle\footnote{You can draw this shape away from the polygon, so it's ok if your line segment between the two dots doesn't cross the side of the triangle anymore!  If it doesn't seem possible:  go back to Step 1, draw your triangles differently, and try again.}.  Now you might have a picture like on the right in Figure \ref{fig:drawing}.  Congratulations! You've drawn your first \emph{tropical curve}\footnote{Unless you've drawn one before.  But hopefully it was fun anyway!}.  

Tropical curves, and more generally \emph{tropical varieties}, are geometric shapes that can be defined by familiar equations called \emph{polynomials}.  However, these polynomials are interpreted using different rules of arithmetic than usual addition and multiplication, replacing addition with taking a maximum and multiplication with addition.  The study of these shapes is called \emph{tropical geometry}, although we can also study other areas of mathematics with these new rules of arithmetic.  In general, we call these subjects \emph{tropical mathematics}.

The first question most people have about tropical mathematics is why it is called ``tropical''.  One of the pioneers of tropical mathematics was Imre Simon, a mathematician and computer scientist who was a professor at the University of S\~{a}o Paulo Brazil.    The adjective \emph{tropical} to describe the field was coined by French mathematicians (Dominique Perrin  or Christian Choffrut, depending on who you ask \cite{pin,simon}) in Professor Simon's honor, based on the proximity of his university to the Tropic of Capricorn.

The second question most people have is why on Earth we would ever redefine our rules of arithmetic in this way.  It turns out that it leads to some incredibly useful and beautiful mathematics.  The first applications of this max-plus arithmetic were in the world of tropical linear algebra, where studying matrix multiplication and related problems in this setting helped solve automation and scheduling problems. More recently, tropical geometry arose as a skeletonized version of algebraic geometry, a major area of mathematics that studies solutions to polynomial equations.  By ``tropicalizing'' solution sets to polynomial equations, we can turn algebro-geometric problems into combinatorial ones, studying more hands-on objects and then lifting that information back to the classical world.  Beyond having applications to computational algebraic geometry, this has allowed for theorems, some new and some old, to be proven in a purely tropical way.

\subsection{Tropical Arithmetic and Tropical Linear Algebra}

The set of real numbers $\mathbb{R}$, equipped with addition $+$ and multiplication $\times$, has the algebraic structure of a \emph{field}.  This means we can add, subtract, multiply, and divide (except by zero), and that arithmetic works essentially how we expect it to. For instance, there's an additive identity $0$, which doesn't change anything when added to it; and there's a multiplicative identity $1$, which doesn't change anything when multiplied by it.  The operations also play well together:  for any $a,b,c\in\mathbb{R}$, we have $a\times (b+c)=a\times b+a\times c$.  If we forget about the fact that we can divide for a minute, all these properties (together with commutativity and associativity of our operations) mean that $(\mathbb{R},+,\times)$ is a \emph{commutative ring with unity}.

Let's now redefine arithmetic on the real numbers with \emph{tropical addition} $\oplus$ and \emph{tropical multiplication} $\odot$, where $a\oplus b=\max\{a,b\}$ and $a\odot b=a+b$.  So, $2\oplus 3=3$ and $2\odot 3=5$.  Instead of only allowing real numbers, we use the slightly larger  set $\overline{\mathbb{R}}=\mathbb{R}\cup\{-\infty\}$, where $-\infty$ has the property that it is smaller than any element of $\mathbb{R}$.  This means, for instance, that $-\infty\oplus 2=2$, and $-\infty \odot 2=-\infty$.

The triple $(\overline{\mathbb{R}},\oplus,\odot)$ \emph{almost} has the structure of a commutative ring with unity, with $-\infty$ as the additive identity and $0$ as the multiplicative identity.  However, elements do not have additive inverses. The equation $1\oplus x=0$ has no solution, since we cannot ``subtract'' $1$ from both sides.  Thus, the triple  $(\overline{\mathbb{R}},\oplus,\odot)$ is a \emph{semiring}, and in particular we call it the {\emph{tropical semiring}}\index{tropical semiring}\footnote{We could have just as easily defined tropical addition as taking the minimum of two numbers.  (Instead of $-\infty$, we would have used $\infty$ as our additive identity.)  Some researchers use the min convention, which is especially useful when studying connections to algebraic geometry; others use the max convention, which is more useful for highlighting certain dualities.  Pay attention to the introductions of books and papers to determine which convention they're using!}.

\begin{exercise}  Verify that tropical addition and tropical multiplication satisfy the law of distributivity.  That is, show that for any $a,b,c\in\overline{\mathbb{R}}$, we have $a\odot(b\oplus c)=(a\odot b)\oplus (a\odot c)$.  Then explain why every element of $\overline{\mathbb{R}}$, besides the additive identity, has a multiplicative inverse.  Because of this it would also be reasonable to refer to $(\overline{\mathbb{R}},\oplus,\odot)$ as the  \emph{tropical semifield}.
\end{exercise}

Historically, the first use of these max-plus operations as an alternative to plus-times came in the world of {\emph{max-linear algebra}}\index{max-linear algebra}, which is similar to linear algebra over the real numbers except that all instances of $+$ and $\times$ are replaced with $\oplus$ and $\odot$.  An example of matrix multiplication with these operations would be

\begin{equation}
\left(\begin{smallmatrix}5&2\\-1&8\end{smallmatrix}\right)\odot \left(\begin{smallmatrix}1&0\\2&-\infty\end{smallmatrix}\right)=\left(\begin{smallmatrix}(5\odot 1)\oplus (2\odot 2)&(5\odot 0)\oplus (2\odot -\infty)\\(-1\odot 1)\oplus(8\odot 2)&(-1\odot 0)\oplus(8\odot -\infty)\end{smallmatrix}\right)=\left(\begin{smallmatrix}6&5\\10&-1\end{smallmatrix}\right).
\end{equation}

There are many natural questions, equations, or definitions coming from usual linear algebra that, when studied tropically, boil down to a scheduling, optimization, or feasibility problem.  We list a few here, and refer the reader to \cite{butkovic} for more details:

\begin{itemize}
\item  Solving equations of the form $A\odot \textbf{x}\leq \textbf{b}$, where $A$ and $\textbf{b}$ are given, solves a scheduling problem.
\item  Finding the determinant of a matrix solves a job assignment problem.  (We have to be careful what we mean by ``determinant'', since there are no negatives tropically!)
\item  Finding an eigenvalue of a matrix finds the shortest weighted cycle on the weighted graph given by the matrix.  (And strangely, this matrix only has that one eigenvalue.)
\end{itemize}

\chproblem{Explain why each of the above linear algebra topics has the  given interpretation when working tropically.}

\resproject{Study the complexity of tropical matrix multiplication.  For both tropical and classical matrix multiplication, the usual algorithm for multiplying two $n\times n$ matrices (namely taking the dot product of rows and columns) uses $n^3$ multiplications.  However, an algorithm for classical matrix multiplication due to Strassen \cite{strassen} has a runtime of $O(n^{2.807})$, with more recent algorithms pushing the runtime down to $O(n^{2.3728639})$ \cite{legall}.  Can such improvements be made for tropical matrix multiplication?

More generally, study the computational complexity of problems in max-linear algebra.}

\subsection{Tropical Polynomials and Tropical Varieties}

A traditional {\emph{polynomial}}\index{polynomial} in $n$ variables over $\mathbb{R}$ is a sum of terms, each of which consists of a coefficient from $\mathbb{R}$ multiplied by some product of those $n$ variables (possibly an empty product; possibly with repeats).  We study the set of points where these polynomials \emph{vanish}; in other words, we set these polynomials equal to $0$, and study the solution sets in $\mathbb{R}^n$.

\begin{example}  The polynomial $x^2-5x+6$, the polynomial $x^2+y^2-1$, and the polynomial $x^2+y^2+z^2-1$ are polynomials in one, two, and three variables, respectively.  The solution sets obtained by setting these polynomials equal to $0$ are the finite set $\{2,3\}$ in $\mathbb{R}$; the unit circle in $\mathbb{R}^2$; and the unit sphere in $\mathbb{R}^3$, respectively.

Note that the solution set of $\{2,3\}$ to $x^2-5x+6=0$ (usually referred to as the \emph{roots} of the polynomial) gives a factorization, namely $x^2-5x+6=(x-2)(x-3)$.  This illustrates the Fundamental Theorem of Algebra:  that any non-constant polynomial in one variable can be factored into linear terms, each of the form $x-\alpha$ with $\alpha$ a root\footnote{There is a bit more fine print:  we must work over $\mathbb{C}$, the field of complex numbers, which is algebraically closed; and we may have to include multiple copies of the same term, based on the \emph{multiplicity} of the root.}.
\end{example}

\emph{Algebraic geometry} is the field of mathematics that studies shapes defined by the vanishing of polynomials. \emph{Tropical geometry}, in parallel, studies shapes defined by \emph{tropical polynomials}.  {Tropical polynomials}\index{tropical polynomial} are the same as usual polynomials, except with all addition and multiplication replaced with tropical addition and tropical multiplication. This includes multiplication of variables, so that $x^2 y$ is interpreted as $x\odot x\odot y=x+x+y=2x+y$.

\begin{figure}[hbt]
\sidecaption
\includegraphics[scale=.65]{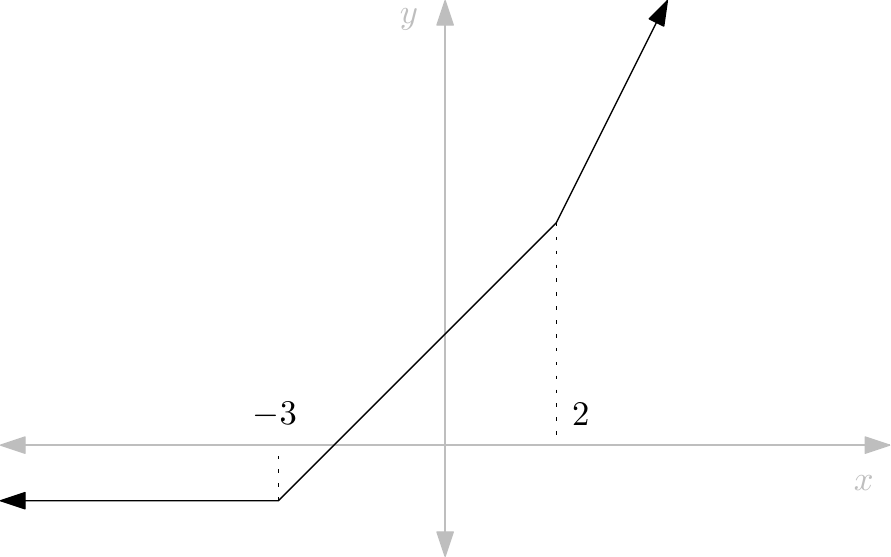}
\caption{The graph of the tropical polynomial $x^2\oplus (2 \odot x) \oplus (-1)$}
\label{fig:one_variable_graph}       
\end{figure}

\begin{example}  The tropical polynomial in one variable $x^2\oplus (2 \odot x) \oplus (-1)$ can be written in classical notation as $\max\{2x,x+2,-1\}$.  The graph of this polynomial, interpreted as a function from $\mathbb{R}$ to $\mathbb{R}$, is illustrated in Figure \ref{fig:one_variable_graph}.
\label{ex:one_variable}
\end{example}

Although we could set a tropical polynomial equal to $0$, the resulting solution set would not be especially meaningful:  most tropical polynomials in one variable are equal to $0$ at at most one point, which doesn't give much information about the polynomial.  Instead, we study the points where the \textbf{maximum is achieved (at least) twice}. In the polynomial from Example \ref{ex:one_variable}, the maximum is achieved twice at two points: when $x=-3$ (where the $2\odot x$ and $-1$ terms tie for the maximum), and when $x=2$  (where the $x^2$ and $2\odot x$ terms tie for the maximum).

\begin{definition} \label{def:vanish} We say that a tropical polynomial $p(x_1,\ldots,x_n)$ {\emph{vanishes}}\index{tropical vanishing} at a point $(a_1,\ldots,a_n)$ if the maximum in $p(a_1,\ldots,a_n)$ is achieved at least twice.  If $p(x)$ is a tropical polynomial in one variable that vanishes at $a$, we say that $a$ is a \emph{root} of $p(x)$.
\end{definition}

As with classical roots, we can give tropical roots a notion of multiplicity:  it is the change in slope going from one linear portion of the graph to the next at that root.  So in Example \ref{ex:one_variable}, both roots have multiplicity $1$, since the slope changes from $0$ to $1$, and then from $1$ to $2$.

\begin{exercise}  We say that a tropical polynomial in one variable has a \emph{root at $-\infty$} if the leftmost linear part of its graph does not have slope $0$; the \emph{multiplicity} of that root is defined to be the slope of that ray.  With this definition, prove that any tropical polynomial in one variable of degree $n$ has exactly $n$ roots in $\overline{\mathbb{R}}$, counted with multiplicity.  (In this sense, $\overline{\mathbb{R}}$ is ``tropically algebraically closed.'')
\end{exercise}

A natural question to ask is whether the tropical roots of a tropical polynomial in one variable have any real meaning.  At least in our example, they give information about how to factor the polynomial:  the reader can verify that $x^2\oplus (2 \odot x) \oplus (-4)=(x\oplus -3)\odot(x\oplus 2)$.   This property holds in general, if we are willing to consider factorizations that give the correct function, even if not the correct polynomial.  (Check and see why $x^2\oplus 0$ and $x^2\oplus (-100\odot x)\oplus 0$ define the same function, even though they're different polynomials!)

\chproblem{Prove the {Tropical Fundamental Theorem of Algebra}\index{The Tropical Fundamental Theorem of Algebra}:  that any tropical polynomial $p(x)$ in one variable is equal, as a function, to
\begin{equation}c\odot (x\oplus \alpha_1)^{\mu_1}\odot(x\oplus \alpha_2)^{\mu_1}\odot \cdots\odot(x\oplus\alpha_k)^{\mu_k},
\end{equation}
where $\alpha_1,\ldots,\alpha_k$ are the tropical roots of $p$, with multiplicities $\mu_1,\ldots,\mu_k$, respectively, and where $c$ is a constant.  
}

\resproject{Study the {factorization of tropical polynomials}\index{tropical polynomial factorization} in more than one variable.  Work in this direction has been done in \cite{lin-tran}, who provide efficient algorithms for certain classes of polynomials, even though in general this is an NP-complete problem.}

Moving beyond polynomials in just one variable, we obtain tropical vanishing sets more complex than finite collections of points.  In Section \ref{section:plane} we study tropical polynomials in two variables in depth, as well as the \emph{tropical curves} they define in $\mathbb{R}^2$.  In Section \ref{section:three_dimensions} we consider tropical polynomials in three variables, which define tropical surfaces.  We also describe how intersecting such surfaces can give rise to tropical curves in three dimensions. In Section \ref{section:tropicalization} we discuss the connection between algebraic geometry and tropical geometry through the tool of tropicalization.
\subsection{Some Tropical Resources}

Throughout this chapter we provide many references to books and articles on tropical geometry, both as sources for results and as great places to find ideas for research projects.  We will frequently reference \emph{An Introduction to Tropical Geometry} by Maclagan and Sturmfels \cite{maclagan-sturmfels}, a graduate text that thoroughly develops the structure of tropical varieties and their connection to algebraic geometry.  That book uses the min convention, while we use the max convention, so we adapt their results as necessary.

The material presented in this chapter, as well as in \cite{maclagan-sturmfels}, looks at tropical geometry from an \emph{embedded} perspective, where tropical varieties are subsets of Euclidean space.  Another fruitful avenue is to look at tropical varieties, especially tropical curves, from an \emph{abstract} perspective, under which tropical curves are thought of as graphs, possibly with lengths assigned to the edges.  In the case of graphs without edge lengths, this theory is thoroughly explored in \cite{chip-firing}.  We also refer the reader to \cite{baker-norine, chan2, cd, gk, mz} for research articles incorporating this perspective.

Finally, there are many fantastic computational tools that help in exploring tropical geometry, both for computing examples and for implementing algorithms.  Here are a few that we'll reference in this chapter, all free to download:

\begin{itemize}
\item  \texttt{Gfan} \cite{gfan}, a software package for computing Gr\"{o}bner fans and tropical varieties.
\item \texttt{Macaulay2} \cite{m2},  a computer algebra system.  Especially useful for us are the \texttt{Polyhedra} and \texttt{Tropical} packages.
\item \texttt{polymake} \cite{polymake}, which is open source software for research in polyhedral geometry. Among many other things, it can deal with polytopes and tropical hypersurfaces.
\item \texttt{TOPCOM} \cite{topcom}, a package for computing Triangulations Of Point Configurations and Oriented Matroids.  As we'll see in Sections \ref{section:plane} and \ref{section:three_dimensions}, being able to find triangulations of polygons and polytopes goes hand in hand with researching tropical varieties.
\end{itemize}

\section{Tropical Curves in the Plane}
\label{section:plane}

Let $p(x,y)$ be a tropical polynomial in two variables with at least two terms.  Let $S$ be the set of all pairs $(i,j)\in\mathbb{Z}^2$ such that a term of the form $c_{ij}\odot x^i\odot y^j$ appears in $p(x,y)$ with $c_{ij}\neq-\infty$; in other words, $S$ is the set of all exponent pairs that actually show up in $p(x,y)$.  We can then write our polynomial as
\begin{equation}
p(x,y)=\bigoplus_{(i,j)\in S}c_{ij}\odot x^i\odot y^j,
\end{equation}
or in classical notation as
\begin{equation}p(x,y)=\max_{(i,j)\in S}\{c_{ij}+ix+jy\}.\end{equation}
As established in Definition \ref{def:vanish}, we say $p(x,y)$ vanishes at a point if this maximum is achieved at least twice at that point.  We call the set of points in $\mathbb{R}^2$ where $p$ vanishes the {\emph{tropical curve}}\index{tropical curve} defined by $p$.  Let $\mathcal{T}(p)$ denote this tropical curve.

\begin{example}

Let $p(x,y)=x\oplus y\oplus 0$. Written in classical notation, $p(x,y)=\max\{x,y,0\}$.  The maximum in this expression is achieved at least twice if two of the terms are equal, and greater than or equal to the third.  This occurs at the point $(0,0)$\footnote{In fact, the maximum occurs three times at this point.}, and along three rays emanating from this point:  when $x=y\geq 0$, when $x=0\geq y$, and when $y=0\geq x$.  The tropical curve $\mathcal{T}(p)$ is illustrated in Figure \ref{fig:line}.  As mentioned in Exercise \ref{exercise:lines}, we call this tropical curve a tropical line.

\begin{figure}[hbt]
\sidecaption
\includegraphics[scale=.65]{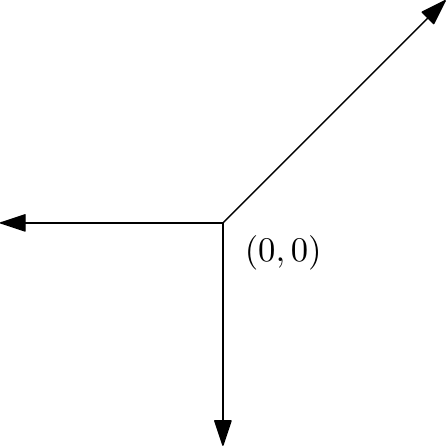}
\caption{The tropical line defined by $x\oplus y\oplus 0$}
\label{fig:line}       
\end{figure}

\end{example}

\begin{exercise}\label{exercise:lines}  Any tropical curve defined by a tropical polynomial of the form $a\odot x\oplus b\odot y\oplus c$, with $a,b,c\in\mathbb{R}$, is called a {\emph{tropical line}}\index{tropical line}.  Determine all the possibilities for what a tropical line can look like.  What if we allow one of $a,b,$ or $c$ to be $-\infty$?
\end{exercise}

\subsection{Convex Hulls and Newton Polygons}

A set in $\mathbb{R}^2$ (or more generally in $\mathbb{R}^n$) is called \emph{convex} if any line segment connecting two points in the set is also contained in the set.  The \emph{convex hull} of a collection of points is the ``smallest'' convex set containing all the points\footnote{More formally, it is the intersection of all convex sets containing the points. See if you can prove that such an intersection is still convex!}.   The {\emph{Newton polygon}}\index{Newton polygon} of $p(x,y)$, written $\textrm{Newt}(p)$, is the convex hull of all the points in $S$.  That is,
\begin{equation}\textrm{Newt}(p)=\textrm{conv}\left(\{(i,j)\in\mathbb{Z}^2\,|\,x^i\odot y^j\textrm{appears in $p(x,y)$ with $c_{ij}\neq-\infty$}\}\right).\end{equation}
As the convex hull of finitely many points in $\mathbb{R}^2$,  $\textrm{Newt}(p)$ is either empty, a point, a line-segment, or a two-dimensional polygon.  To avoid certain trivial cases, we'll assume that we've chosen $p$ such that $\textrm{Newt}(p)$ is a two-dimensional polygon.  It is a {\emph{lattice polygon}}\index{lattice polygon}, meaning that all vertices are \emph{lattice points}, which are are points with integer coordinates.  In the special case that $\textrm{Newt}(p)=\textrm{conv}\{(0,0),(d,0),(0,d)\}$ for some positive integer $d$, we say that the polynomial has degree $d$, and we call the Newton polygon the \emph{triangle of degree $d$}, denoted $T_d$.

\begin{figure}[hbt]
\sidecaption
\includegraphics[scale=.65]{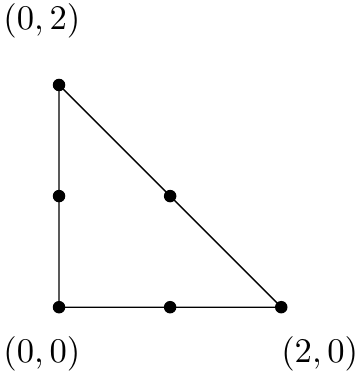}
\includegraphics[scale=.65]{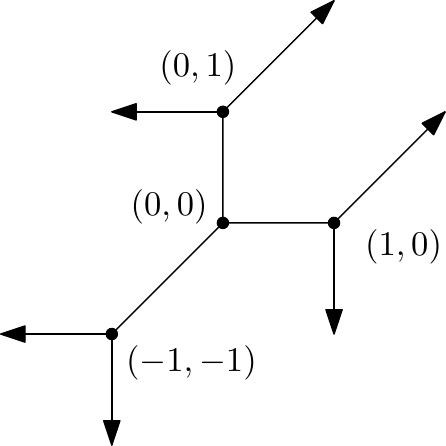}
\caption{The Newton polygon of $(1\odot x^2)\oplus (1\odot y^2)\oplus (2\odot x y)\oplus (2\odot x)\oplus (2\odot y)\oplus 1$, along with the tropical curve the polynomial defines}
\label{fig:quadric_polygon_and_curve}       
\end{figure}

\begin{example}  Let $p(x,y)=(1\odot x^2)\oplus (1\odot y^2)\oplus (2\odot x y)\oplus (2\odot x)\oplus (2\odot y)\oplus 1$.  Then we have that $S=\{(2,0),(0,2),(1,1),(1,0),(0,1),(0,0)\}$, so $\textrm{Newt}(p)$ is the triangle of degree $2$, and $p(x,y)$ is a polynomial of degree $2$. The Newton polygon, along with the tropical curve $\mathcal{T}(p)$, are illustrated in Figure \ref{fig:quadric_polygon_and_curve}.   Some preliminary connections between $\textrm{Newt}(p)$ and $\mathcal{T}(p)$ can already be observed:  the rays in $\mathcal{T}(p)$ point in directions that are perpendicular and outward relative to the edges of $\textrm{Newt}(p)$.  However, there are other features of the tropical curve not visible from the Newton polygon; for instance, there are three bounded edges, and there are four \emph{vertices}, where multiple edges or rays come together.

\label{ex:quadratic}
\end{example}

\subsection{Subdivisions and the Duality Theorem}

Since it was presented without justification, the reader might wonder:  how did we determine  $\mathcal{T}(p)$  in Example \ref{ex:quadratic}?  One brute force way could be to take every possible pair among the $6$ terms in $p(x,y)$ (there are $15$ ways to do this), set them equal to each other, and try to determine whether those two terms ever tie for the maximum, and if so, where.  It turns out that studying the Newton polygon of $p$ leads to a much more elegant approach.

Let $P$ be a lattice polygon, and $S=P\cap \mathbb{Z}^2$ be the set of integer coordinate points in $P$.  Let $h:S\rightarrow\mathbb{R}$ be any function assigning real number values\footnote{This definition will still work even if we define $h:S\rightarrow\mathbb{R}\cup\{-\infty\}$, as long as $h$ does not map any vertices of $P$ to $-\infty$.} to each element of $S$; we refer to $h$ as a \emph{height function}.  
We then define a set $A$ of points in $\mathbb{R}^3$ by ``lifting'' the points of $S$ to the heights prescribed by $h$:
\begin{equation}A=\{(i,j,h(i,j))\,|\, (i,j)\in S\}.\end{equation}
Take the convex hull of $A$ in $\mathbb{R}^3$. Unless all the points of $A$ lie on a plane, this convex hull is a three-dimensional \emph{polytope}, the three-dimensional analog of a polygon, whose boundary consists of two-dimensional  polygonal faces meeting along edges.  Viewed from above, $\textrm{conv}(A)$ looks like $P$, except subdivided by these upper polygonal faces.  We call this subdivision of $P$ the {\emph{subdivision induced by $h$}}\index{induced subdivision}.  The faces of $\textrm{conv}(A)$ that are visible from above form the {\emph{upper convex hull}}\index{upper convex hull} of $A$.

\begin{example}  Let $p(x,y)$ be as in Example \ref{ex:quadratic}.  Let $P=\textrm{Newt}(p)$, and $S=P\cap\mathbb{Z}^2$.  Define $h:S\rightarrow\mathbb{R}$ using the coefficients of $p(x,y)$, so that $h(i,j)=c_{i,j}$.  Then the set $A$ consists of the six points $\{(0,0,1),(1,0,2),(2,0,1),(0,1,2),(1,1,2),(0,2,1)\}$, illustrated on the left in Figure \ref{fig:convex_hull}.  Their convex hull is then a polytope with $8$ triangular faces, illustrated in the middle of the figure.  Of these faces, the $4$ that are colored are visible from above, giving the induced subdivision of $P$shown towards the right.  The tropical curve $\mathcal{T}(p)$ is reproduced, with vertices colored the same as their corresponding triangles, as described in Theorem \ref{theorem:duality} below.

\begin{figure}[hbt]
\sidecaption
\includegraphics[scale=.65]{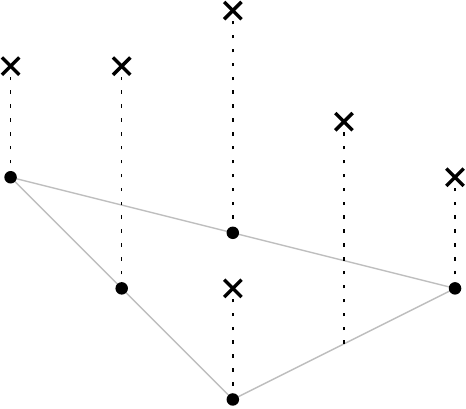}\qquad
\includegraphics[scale=.65]{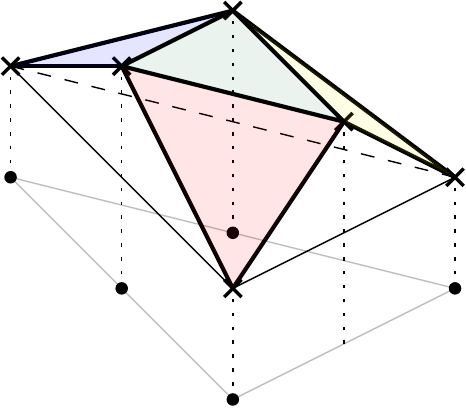}\qquad
\includegraphics[scale=.65]{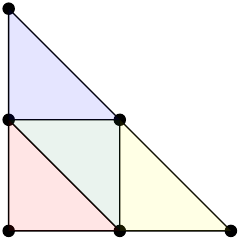}
\includegraphics[scale=.5]{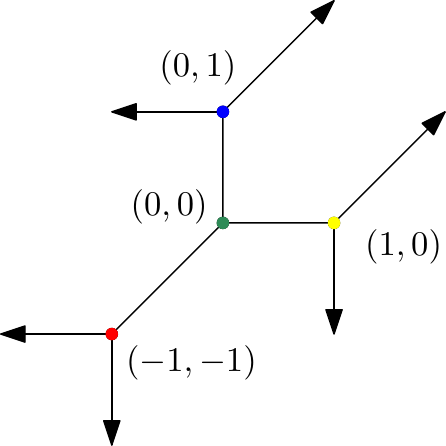}
\caption{The  points of $A$ labelled as $\times$'s, their convex hull, the induced subdivision of the triangle, and the dual tropical curve}
\label{fig:convex_hull}       
\end{figure}

\label{ex:quadratic2}
\end{example}

The subdivision of the Newton polygon induced by the coefficients of the tropical polynomial gives us almost all the information regarding how to draw the tropical curve in the plane.  Although this result holds in much more generality, we spell it out explicitly in the case of two-variables.

\begin{theorem}[{The Duality Theorem}\index{The Duality Theorem}, {\cite[Proposition 3.1.6]{maclagan-sturmfels}}]  Let $p(x,y)$ be a tropical polynomial with $P=\textrm{Newt}(p)$ two-dimensional.  Then the tropical curve $\mathcal{T}(p)$ is \emph{dual} to the subdivision of $P$ induced by the coefficients of $p(x,y)$ in the following sense:

\begin{itemize}
\item  Vertices of $\mathcal{T}(p)$ correspond to polygons in the subdivision of $P$.
\item  Edges of $\mathcal{T}(p)$ correspond to interior edges in the subdivision of $P$.
\item  Rays of $\mathcal{T}(p)$ correspond to boundary edges in the subdivision of $P$.
\item  Regions of $\mathbb{R}^2$ separated by $\mathcal{T}(p)$ correspond to lattice points of $P$ used in the subdivision.
\end{itemize}
Moreover, two vertices of $\mathcal{T}(p)$ are connected by an edge if and only if their corresponding polygons in the subdivision share an edge, and the edge in the Newton polygon is perpendicular to the edge in the subdivision; and the rays emanating from a vertex in $\mathcal{T}(p)$ correspond to boundary edges of the corresponding polygon in the subdivision, with the rays in the outward perpendicular directions to the boundary edges of $P$.

\label{theorem:duality}
\end{theorem}

So once we have found the subdivision of our Newton polygon, we know exactly what the tropical curve will look like, up to scaling edge lengths and up to translation.  If we find the subdivision from Example \ref{ex:quadratic2}, then our tropical curve could be either of the ones illustrated in Figure \ref{fig:two_curves} (or infinitely many others!).  However, we can nail down the exact coordinates of the vertices by solving for the relevant three-way-ties.  For instance, the top-most vertex of the tropical curve corresponds to the triangle with vertices at $(0,2)$, $(0,1)$, and $(1,1)$ in the subdivision, so the coordinates of the vertex are located at the (unique) three-way tie between the $y^2$, the $y$, and the $x y$ terms.

\begin{figure}[hbt]
\sidecaption
\includegraphics[scale=.65]{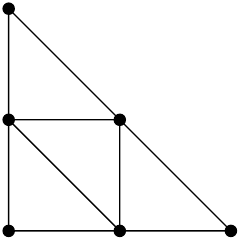}
\includegraphics[scale=.65]{quadratic_curve}
\includegraphics[scale=.65]{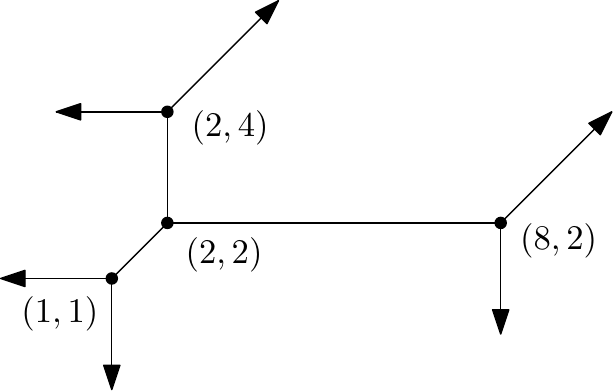}
\caption{A subdivision of a Newton polygon, and two possible tropical curves dual to it}
\label{fig:two_curves}       
\end{figure}

\begin{figure}[hbt]
\sidecaption
\includegraphics[scale=.8]{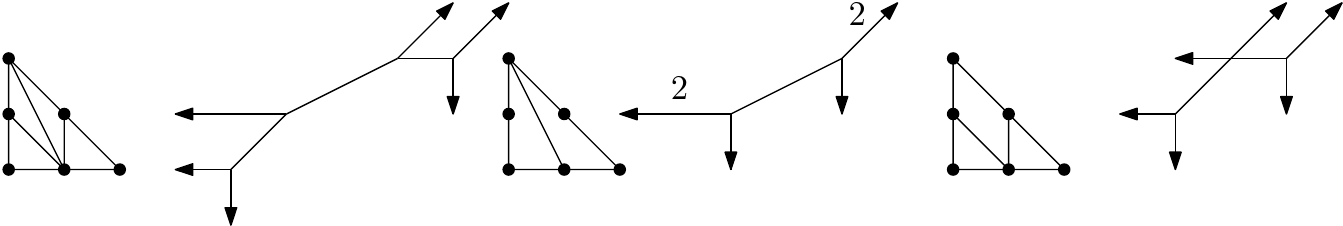}
\caption{Three tropical curves with the same Newton polygon, dual to different subdivisions.  The first tropical curve is smooth, and the other two are not.}
\label{fig:three_curves}       
\end{figure}

Sometimes there is information present in the polynomial or in the subdivision of the Newton polygon that is lost in the tropical curve.  For instance, if $p(x,y)=x^2\oplus y^2\oplus 0$, then $\mathcal{T}(p)$ is, as a set, the tropical line from Figure \ref{fig:line}.  By only considering this tropical curve as a set, we thus lose information about the starting polynomial.  This leads us to decorate the edges and rays of our tropical curves with \emph{weights}.  In particular, each edge or ray is given a positive integer weight $m$, where $m$ is equal to one less than the number of lattice points on the dual edge of the subdivision.  Several tropical curves with the same Newton polygon are illustrated in Figure \ref{fig:three_curves}, with all weights above $1$ labelled.  If a tropical curve has all weights equal to $1$, and each vertex has a total of three edges and rays emanating from it, then we call the {tropical curve \emph{smooth}}\index{smooth tropical curve}.  Equivalently, a tropical curve is smooth if its dual subdivision is a {\emph{unimodular triangulation}}\index{unimodular triangulation}, meaning that every polygon in the subdivision is a triangle with no lattice points besides its vertices\footnote{By Pick's Theorem \cite{pick}, this in turn is equivalent to every polygon in the subdivision being a triangle with area $1/2$.}.

\begin{exercise}  Let $p(x,y)$ be a tropical polynomial of degree $d$ such that  $\mathcal{T}(p)$ is smooth.  Determine the number of edges, rays, and vertices of $\mathcal{T}(p)$.  (\textbf{Hint}:  count up the corresponding objects in a unimodular triangulation of the triangle $T_d$.  You can use the fact that any triangle in such a triangulation has area $1/2$.)
\end{exercise}

\chproblem{Show that any tropical curve satisfies the following {\emph{balancing condition}}\index{balancing condition}\footnote{This is a special case of a much more general result called the \emph{Structure Theorem}, which says that any tropical variety has the structure of a weighted, balanced polyhedral fan of pure dimension.  See \cite[Theorem 3.3.5]{maclagan-sturmfels}.}:  choose a vertex, and let $\left<a_1,b_1\right>,\left<a_2,b_2\right>,\ldots ,\left<a_\ell,b_\ell\right>$ be the outgoing directions of the rays and edges emanating from the vertex, where $a_i,b_i\in\mathbb{Z}$ and $\gcd(a_i,b_i)=1$ for all $i$.  Let $m_i$ denote the weight of the $i^{th}$ edge/ray.  Show that $m_1\times \left<a_1,b_1\right>+m_2\times \left<a_2,b_2\right>+\cdots +m_\ell\times\left<a_\ell,b_\ell\right>=\left<0,0\right>$.}

\begin{exercise}
Consider the subset $C$ of $\mathbb{R}^2$ illustrated in Figure \ref{fig:show_tropical_curve}. It consists of three rays, all emanating from the origin, in the directions $\left<1,0\right>$, $\left<0,1\right>$, and $\left<-2,-1\right>$.  Show that $C$ is a tropical curve by finding a tropical polynomial $p(x,y)$ such that $C=\mathcal{T}(p)$.  (\textbf{Hint}:  the previous Challenge Problem might be useful!)

\begin{figure}[hbt]
\sidecaption
\includegraphics[scale=.7]{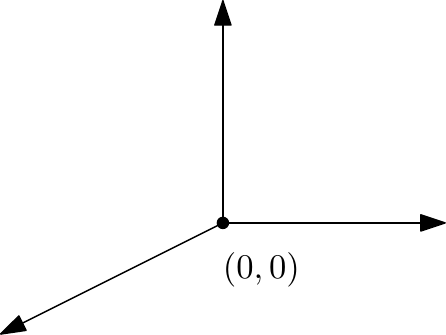}
\caption{A set that turns out to be a tropical curve}
\label{fig:show_tropical_curve}       
\end{figure}
\end{exercise}

Armed with our Duality Theorem, one way to study tropical curves is the following:  choose a tropical polynomial, find the induced subdivision of its Newton polygon, and draw it, solving for the exact coordinates of the vertices.  Perhaps the most challenging step is finding the induced subdivision; this can be accomplished with such computational tools as \texttt{polymake}, \texttt{TOPCOM}, and \texttt{Macaulay2}.  

\begin{figure}[hbt]
\sidecaption
\includegraphics[scale=.65]{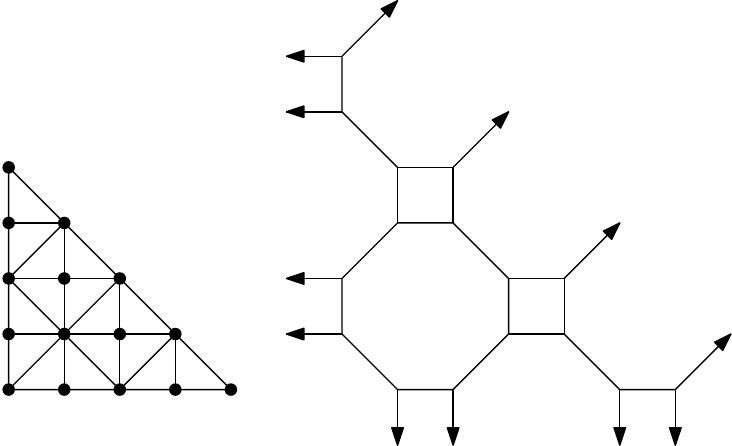}
\qquad
\includegraphics[scale=.65]{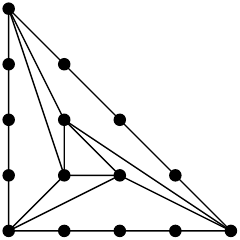}
\caption{A regular subdivision with a dual tropical curve, and a non-regular subdivision}
\label{fig:genus_3}       
\end{figure}

Here we take another approach, similar to the very start of this chapter.  Rather than starting with a tropical polynomial, choose the Newton polygon, and simply draw a subdivision, perhaps a unimodular triangulation.  Then try to draw a tropical curve dual to it. (This is exactly the method from the start of Section 1.) An example of a triangulation of the triangle of degree $4$ is illustrated in Figure \ref{fig:genus_3}, along with a tropical curve that is dual to it.  Note that to draw this tropical curve, we never needed to find a tropical polynomial defining it!

Sadly, this approach does not always work.  A tropical curve can be drawn dual to a subdivision if and only if the subdivision is {\emph{regular}}\index{regular subdivision}, meaning that it is induced by some height function.

\exercise{Consider the subdivision on the right in Figure \ref{fig:genus_3}.  Show that it is not a regular triangulation.  (You might argue that no height function could have induced that triangulation; or you could argue that it's impossible to draw a tropical curve dual to it.)}

It turns out that there are $1279$ unimodular triangulations of the triangle of degree $4$ up to symmetry \cite{blmpr, bjms}, and only one of them is non-regular:  it is the unique unimodular triangulation that completes the non-regular subdivision from Figure \ref{fig:genus_3}.  Similar phenomena occur for ``small'' polygons, whereby most triangulations end up being regular, so that drawing dual tropical curves is usually possible.  For larger polygons, regular subdivisions seem to become rarer and rarer. See \cite{kz} for many results in the case that the polygon is a lattice rectangle, as well as \cite{triangulations} for results in a more general setting.

\chproblem{Let $n$ be a positive integer, and let $P$ be a $1\times n$ lattice rectangle.  Prove that any subdivision of $P$ is regular.  How many unimodular triangulations are there of $P$?}

\resproject{Study the number of unimodular triangulations of families of lattice polygons, as was done for lattice rectangles in \cite{kz}.  This can involve finding upper and lower bounds that improve those in the literature.  Study the proportion these unimodular triangulations that are regular.  For all these endeavors, \texttt{polymake} and \texttt{TOPCOM} are fantastically useful computational tools.}

\subsection{The Geometry of Tropical Plane Curves}

Many theorems about classical plane curves have analogs within the tropical world.  A prime example of this is B\'{e}zout's Theorem.

\begin{theorem}[{B\'{e}zout's Theorem}\index{B\'{e}zout's Theorem}]  Let $C$ and $D$ be two smooth algebraic plane curves of degrees $d$ and $e$.  If $C$ and $D$ have no common components, then $C\cap D$ has at most $d\times e$ points.  If we are working in projective space over an algebraically closed field, and counting intersection points with multiplicity, then $C \cap D$ has exactly $d\times e$ points. 
\end{theorem}

As shown in \cite{firststeps}, the same result holds for tropical plane curves, once we determine how to count intersection points with multiplicity, and how to deal with tropical curves that intersect ``badly''.

\begin{definition} \label{def:multiplicity}  Suppose two tropical plane curves $C_1$ and $C_2$ intersect at an isolated point $(a,b)$ that is not a vertex of either curve.  Such a point is called a \emph{transversal intersection}.  Let $\left<u_1,v_1\right>$ and $\left<u_2,v_2\right>$ be integer vectors describing the slopes of the edges or rays of $C$ and $D$ containing $(a,b)$, where $\gcd(u_1,v_1)=\gcd(u_2,v_2)=1$, and let the weights of the edges or rays be $m_1$ and $m_2$.  Then the {\emph{multiplicity}}\index{intersection multiplicity} of $(a,b)$ is
\begin{equation}\mu(a,b):=m_1\times m_2\times\left|\det\left(\begin{smallmatrix}u_1&v_1\\u_2&v_2\end{smallmatrix}\right)\right|.\end{equation}
\end{definition}

\begin{example}  Consider the tropical polynomials 
\begin{equation}f=(-1\odot x^2)\oplus  (x y)\oplus (-1\odot y^2)\oplus x\oplus y\oplus (-1) \end{equation} and 
\begin{equation}g=\left(-\frac{1}{2}\odot x^2\right)\oplus  (1\odot x y)\oplus (-2\odot y^2)\oplus x\oplus y\oplus 0 .\end{equation} They both have the triangle of degree $2$ as their Newton polygon, and have induced subdivisions as illustrated on the left in Figure \ref{fig:two_curves}.
As shown on the right, the tropical curves $\mathcal{T}(f)$ and $\mathcal{T}(g)$ intersect in three points.  The multiplicities of these points can be computed as $1$, $1$, and $2$.

\begin{figure}[hbt]
\sidecaption
\includegraphics[scale=.65]{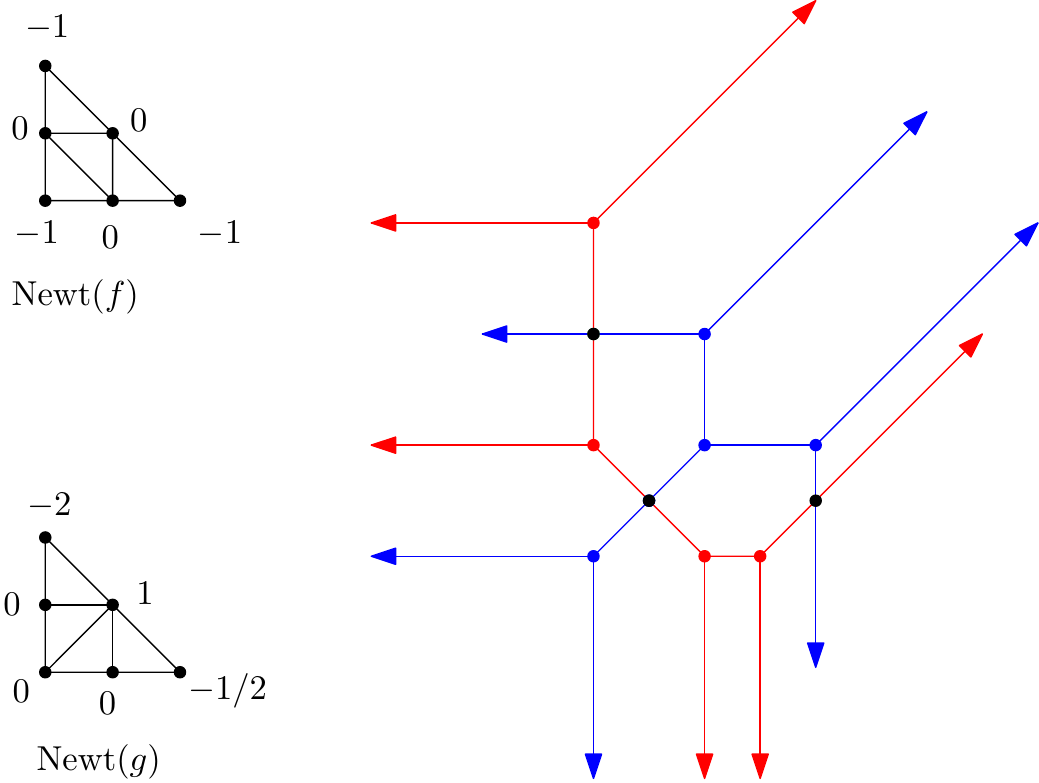}
\caption{The subdivisions induced by $f$ and $g$, and the two tropical curves.}
\label{fig:two_curves}       
\end{figure}

Let's push this example a little further.  If we think of $\mathcal{T}(f)\cup\mathcal{T}(g)$ as $\mathcal{T}(f\odot g)$, then we can consider the dual subdivision of $\textrm{Newt}(f\odot g)$, illustrated in Figure \ref{fig:big_subdivision}.  Every polygon in this subdivision is dual to a vertex of $\mathcal{T}(f\odot g)$, and each vertex in $\mathcal{T}(f\odot g)$ is either  a vertex of $\mathcal{T}(f)$, a vertex of $\mathcal{T}( g)$, or an intersection point.  Note that each polygon dual to an intersection point $(a,b)$ has area equal to $\mu(a,b)$.  

\begin{figure}[hbt]
\sidecaption
\includegraphics[scale=1.52]{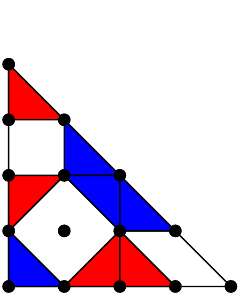}
\caption{The subdivisions induced by $f\odot g$, with blue triangles coming from vertices in $\mathcal{T}(f)$ and red triangles coming from $\mathcal{T}(g)$.}
\label{fig:big_subdivision}       
\end{figure}
\end{example}

\begin{exercise}  Show that if $f$ and $g$ are tropical polynomials of degrees $d$ and $e$, then $f\odot g$ is a tropical polynomial of degree $d\odot e$, and that $\mathcal{T}(f\odot g)=\mathcal{T}(f)\cup\mathcal{T}(g)$.  Then show that the multiplicity of a transversal intersection point of $f$ and $g$ is equal to the area of the corresponding polygon in the subdivision of $\textrm{Newt}(f\odot g)$ induced by $f\odot g$.  
\end{exercise}

\begin{theorem}[{Tropical B\'{e}zout's Theorem}\index{Tropical B\'{e}zout's Theorem}, Transversal Case]  Let $C$ and $D$ be two tropical plane curves of degrees $d$ and $e$ with finitely many intersection points $(a_1,b_1),\cdots,(a_n,b_n)$, all of which are transversal.  Then
\begin{equation}\sum_{i=1}^n\mu(a_i,b_i).\end{equation}
\end{theorem}

Note that we did not need to assume $C$ and $D$ were smooth.  For an even more general result, we need to deal with the possibility that $C$ and $D$ have intersections that are not transversal.  For two tropical curves $C$ and $D$, we compute the {\emph{stable tropical intersection}}\index{stable tropical intersection} as follows.  Let $\textbf{v}=\left<v_1,v_2\right>$ be a vector not parallel to any edge or ray of $C$ and $D$, and for $\varepsilon\in\mathbb{R}^+$ let $D_\varepsilon$ be a translation of $D$ by $\varepsilon\textbf{v}$.  We then define
\begin{equation}C\cap_{st}D=\lim_{\varepsilon\rightarrow 0}C\cap D_{\varepsilon}.\end{equation}
The \emph{multiplicity} of a point in $C\cap_{st}D$ is the sum of the multiplicities of the corresponding points in a small enough perturbation $C\cap D_{\varepsilon}$.

\begin{example}  If $f(x,y)=x\oplus y\oplus 0$ and $g(x,y)=(1\odot x)\oplus y \oplus 0$, then $C=\mathcal{T}(f)$ and $D=\mathcal{T}(g)$ are the tropical lines pictured in Figure \ref{fig:stable_intersection}.  Their set-theoretic intersection is a ray emanating from the point $(-1,0)$.  To find $C\cap_{st}D$, we move $D$ slightly to $D_\varepsilon$, and then move it back to $D$.  In the limit, we find a single stable intersection point at $(-1,0)$.

\begin{figure}[hbt]
\sidecaption
\includegraphics[scale=0.7]{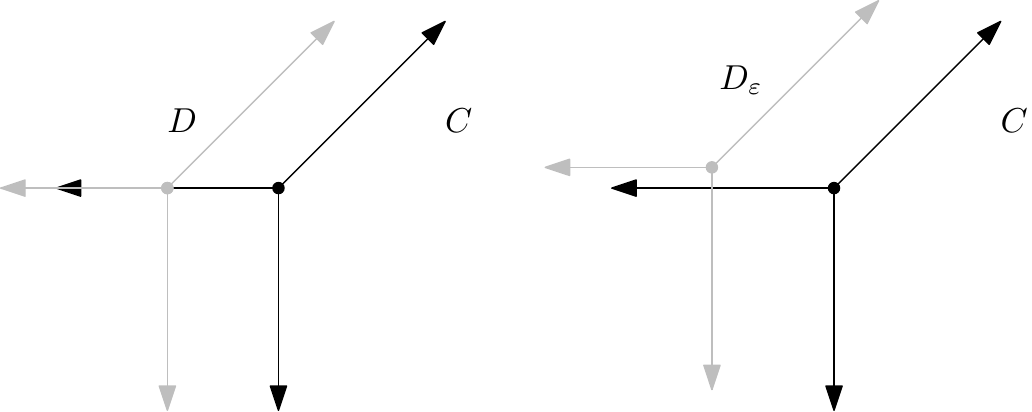}
\caption{Two tropical lines intersecting non-transversally, and a small perturbation used to compute the stable intersection.}
\label{fig:stable_intersection}       
\end{figure}
\end{example}

\begin{exercise}  Show that $C\cap_{st}D$ is a well-defined set of finitely many points, and is independent of the choice of $v$.  Also show that the multiplicity of each point is well-defined.
\end{exercise}

\begin{theorem}[Tropical B\'{e}zout's Theorem, General Case]  Let $C$ and $D$ be two tropical plane curves of degrees $d$ and $e$ with $C\cap_{st}D=(a_1,b_1),\cdots,(a_n,b_n)$.  Then
\begin{equation}\sum_{i=1}^n\mu(a_i,b_i).\end{equation}
\end{theorem}

\chproblem{Prove the transversal case of  tropical B\'{e}zout's Theorem using an area-based argument involving the Newton polygon of $f\odot g$.  Then use this result to prove the general case of tropical B\'{e}zout's Theorem.}

Many classical results about algebraic plane curves involve when two curves are \emph{tangent} to one another at some collection of points.  Recently much work has been done to build up machinery to pose and study these sorts of results in the tropical world.

\begin{definition}  Let $C$ and $D$ be tropical curves.  A \emph{tangency} between $C$ and $D$ is a component of $C\cap D$ such that the stable intersection $C\cap_{st} D$ has more than one point in that component, counted with multiplicity.  We say $C$ and $D$ are \emph{tangent} at that component of $C\cap D$.
\end{definition} 

\begin{figure}[hbt]
\sidecaption
\includegraphics[scale=.75]{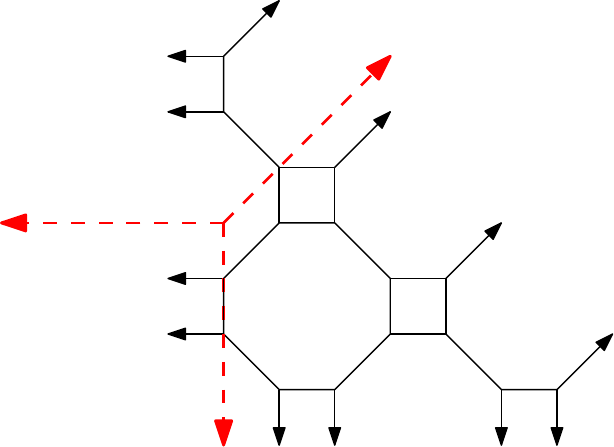}
\caption{A tropical line that is tangent to a tropical curve at two components}
\label{fig:bitangent}       
\end{figure}

A tropical line that is tangent to a degree $4$ curve at two distinct components are illustrated in Figure \ref{fig:bitangent}.  Such an intersection is called a {\emph{bitangent line}}\index{tropical bitangent}, which is also used to refer to an intersection component of multiplicity $4$ or more.

\begin{exercise}  Find all the bitangent lines of the curve from Figure \ref{fig:bitangent}. (Hint:  there are infinitely many of them, but they still admit a nice classification.)
\end{exercise}

Counting bitangent lines is a very classical problem in algebraic geometry.  In 1834, Pl\"{u}cker proved that a smooth algebraic  plane curve of degree $4$  has $28$ bitangent lines \cite{plucker}.  A tropical analog of this fact was proved in \cite{blmpr}.

\begin{theorem}[Theorem 3.9 in \cite{blmpr}]  Let $C$ be a smooth tropical plane curve of degree $4$.  Then $C$ has exactly seven classes\footnote{Loosely speaking, we say two bitangent lines intersecting at $(P,Q)$ and $(P',Q')$ with multiplicity $2$ at each point are equivalent if $(P,Q)$ and $(P',Q')$ are equivalent in the language of \emph{divisor theory}~\cite{gk}.} of bitangent lines.
\end{theorem}

Later work was done to relate this theorem to Pl\"{u}cker's count, starting in \cite{k4} and culminating in \cite{lifting}, which showed how to recover the classical count of $28$ bitangent lines from the tropical count, at least in sufficiently general cases.

\resproject{One great starting point for asking tropical questions is to study tropical versions of algebraic results.  Study, prove, or disprove tropical analogs of these classical results.  You may have to assume something about positions being sufficiently general.

\begin{itemize}
\item  The De Bruijn-Erd\"{o}s theorem \cite{de}:  for any $n$ points not all on a line determining $t$ points,  then $t\geq n$ and if $t=n$, any two lines have exactly one of the $n$ points in common.  (In this latter case, $n-1$ of the points are collinear.)
\item  Steiner's conic problem \cite{bkt}:  given $5$ curves of degree $2$, how many curves of degree $2$ are tangent to all of them? (Classically, the answer is $3264$, although Steiner incorrectly computed it as $7776$.)
\item  The Three Conics Theorem \cite{sevencircles}:  given three conics that pass through two given points, the three lines joining the other two intersections of each pair of conics all intersect at a point.  Dually:  given three conics that share two common tangents,  the remaining pairs of common tangents intersect at three points that are collinear.
\item  The Four Conics Theorem \cite{sevencircles}:   Suppose we are given three conics, where two intersections of each pair lie on a fourth conic.  Then the three lines joining the other two intersections of each pair of conics intersect in a point.
\end{itemize}}

It's also worth determining when tropical geometry does \emph{not} nicely mirror classical algebraic geometry.  We say that an algebraic or a tropical curve $C$ is \emph{irreducible} if it cannot be written as $C_1\cup C_2$, where $C_1\subsetneq C$ and $C_2\subsetneq C$ are curves as well.  One nice property of algebraic curves (and more generally algebraic varieties) is that they admit a unique {decomposition into irreducible components}\index{decomposition into irreducible components} \cite[Theorem 4.6.4]{clo}, just as any integer $n\geq 2$ can be written as a product of primes uniquely (up to reordering).  Tropical curves, however, do not.

\begin{example}\label{example:two_unions}
Consider the set $C$ in $\mathbb{R}^2$ consisting of the (usual) lines $x=0$, $y=0$, and $x=y$.  We claim that $C$ is a tropical curve; you will show this in Exercise \ref{exercise:two_unions} .  We can also write $C$ as $\mathcal{T}(x\oplus y\oplus 0)\cup \mathcal{T}((x y)\oplus x\oplus y)$, or as $\mathcal{T}(x\oplus y)\cup \mathcal{T}(x\oplus 0)\cup \mathcal{T}(y\oplus 0)$, as illustrated in Figure \ref{fig:two_unions}.

\begin{figure}[hbt]
\sidecaption
\includegraphics[scale=0.6]{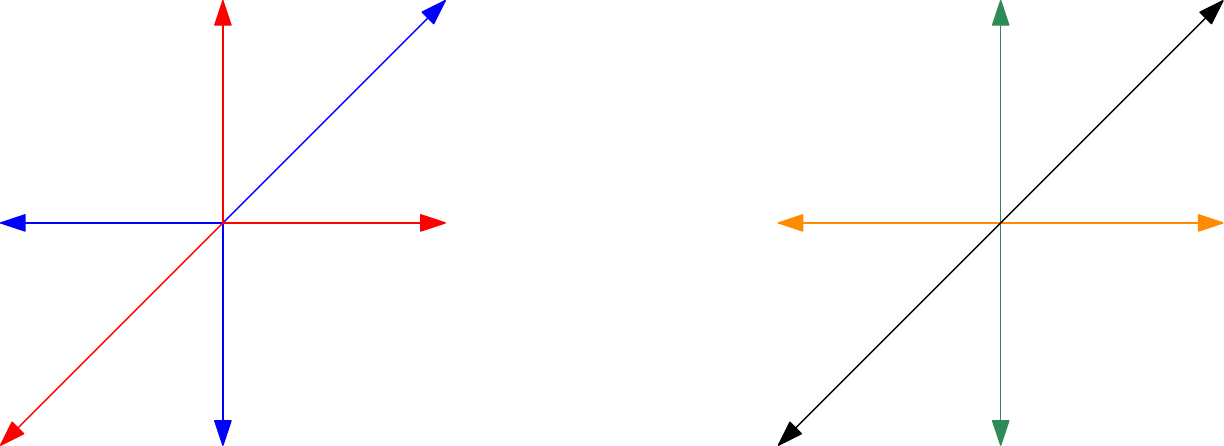}
\caption{A tropical curve that can be decomposed into irreducible tropical curves in two distinct ways.}
\label{fig:two_unions}       
\end{figure}

\end{example}

\begin{exercise} \label{exercise:two_unions} Find a polynomial $f$ such that $C=\mathcal{T}(f)$, where $C$ is the set from Example \ref{example:two_unions}.   How does this polynomial relate to the polynomials defining the two decompositions of $C$ as a union of tropical curves?
\end{exercise}

\resproject{Study how many decompositions a tropical curve can have as a union of tropical curves properly contained within it.  You could stratify this study by the Newton polygon of the curve.  (This is closely related to the research project on factoring tropical polynomials; see if you can see why, especially after you try Exercise \ref{exercise:two_unions}!)}

A new approach in tropical geometry that avoids non-uniqueness of decompositions is to develop \emph{tropical schemes} \cite{tropicalschemes, tropicalideals}, just as algebraic geometers study \emph{algebraic schemes} \cite{hartshorne}.  This model does not consider the tropical curves from the second decomposition in Example \ref{example:two_unions} to be tropical curves, and in fact gives us a unique decomposition in general.

\subsection{Skeletons of Tropical Plane Curves}

\begin{figure}[hbt]
\sidecaption
\includegraphics[scale=.75]{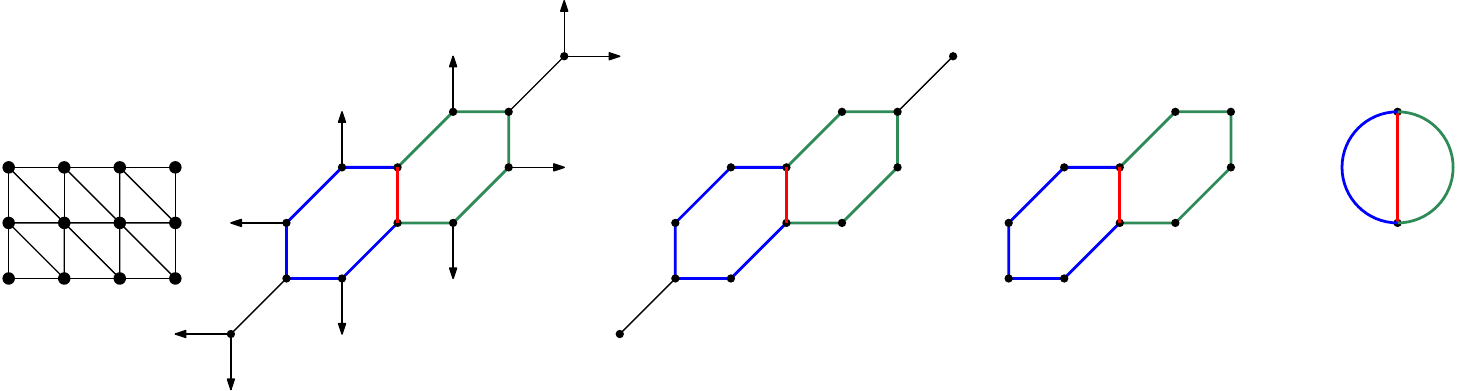}
\caption{A tropical curve (with its dual subdivision) undergoing the process of skeletonization.  The edges of the tropical curve that end up contributing to the skeleton are color-coded based on which final edge they become a part of.}
\label{fig:skeletonization}       
\end{figure}

Choose a lattice polygon $P$ with $g$ interior lattice points, where $g$ is at least $2$.  Write $P_{int}$ for the convex hull of the $g$ interior lattice points; this is either a line segment, or a polygon.  Let $p(x,y)$ be a tropical polynomial with Newton polygon $P$.  Rather than study the full tropical curve $\mathcal{T}(p)$, we can focus on a portion of it called its {\emph{skeleton}}\index{skeleton}.  To find the skeleton, we delete all rays from our tropical curve, and then successively remove any vertices incident to exactly one edge, along with such edges.  This will lead to a collection of vertices and edges, where each vertex is incident to at least two edges.  We ``smooth over'' the vertices incident to two edges, removing such vertices and fusing the two edges into one.  The resulting collection of edges and vertices is called the \emph{skeleton} of the tropical curve.  This process is illustrated in Figure \ref{fig:skeletonization}.

The structure that remains after ``skeletonizing'' a tropical curve is a \emph{graph}\footnote{In fact, there is a bit more structure:  it is a \emph{metric graph}, meaning the edges have lengths.  We'll come back to that later in this subsection.}.  A graph is simply a collection of vertices collected by edges; in our setting, two vertices may be connected to each other by multiple edges, and a vertex may be connected to itself by an edge, which we call a loop. This leads us to the following major question:  Which graphs can appear as the skeleton of tropical plane curve?  To simplify, lets assume that our tropical curves are smooth.

\begin{definition} A graph that is the skeleton of some smooth tropical plane curve is called \emph{tropically planar}, or {\emph{troplanar}}\index{troplanar graph} for short.  The \emph{genus}\footnote{There is another, unrelated definition of \emph{genus} in graph theory, dealing with the smallest number of holes a surface must have to allow a given graph to be embedded on it.} of the graph is the number of bounded regions in the plane formed by a drawing on the graph.  By Euler's formula relating the number of vertices, edges, and faces of a planar graph, we could also define the genus as $E-V+1$ for a graph with $E$ edges and $V$ vertices.
\end{definition}

With these definitions, we can say that the graph on the right in Figure \ref{fig:skeletonization} is troplanar, and has genus $2$.

\begin{exercise} \label{exercise:skeleton_properties} Let $G$ be a troplanar graph.  Show that $G$ is connected (all one piece), planar (able to be drawn in the plane without any edges crossing), and trivalent (meaning that every vertex has three edges coming from it, where a loop counts as two edges).  Also show that the genus of the graph is equal to $g$, the number of interior lattice points of the Newton polygon of any smooth tropical curve that has $G$ as its skeleton.
\end{exercise}

A daunting task is to try to determine which graphs are tropically planar.  Even for fixed $g$, it is not immediately obvious that there is an algorithmic way to do this.  There are several things working in our favor:

\begin{enumerate}
\item There are only finitely many polygons with $g\geq 1$ interior lattice points, up to equivalence\footnote{Here we say two lattice polygons are \emph{equivalent} if one is the image of the other under a matrix transformation $\left(\begin{smallmatrix}a&b\\c&d\end{smallmatrix}\right)$, where $ad-bc=\pm 1$.}.  As discussed in \cite[Proposition 2.3]{bjms}, this follows from results in \cite{scott} and \cite{lz}. An algorithm for finding all such polygons for a given $g$ is presented in \cite{movingout}.
\item  If $P$ and $Q$ are lattice polygons with $P\subset Q$ and $P_{int}=Q_{int}$, all the troplanar graphs arising from $P$ also arise from $Q$ \cite[Lemma 2.6]{bjms}.
\end{enumerate}

\begin{exercise}  Prove item 2 above.
\end{exercise}

Item 1 means that we only need to consider a finite collection of possible Newton polygons for each genus $g$; item 2 decreases that number considerably.  It means that we need to only consider {\emph{maximal} polygons}\index{maximal polygon}, which are those that are not properly contained in any polygon with the same interior lattice points.

\begin{figure}[hbt]
\sidecaption
\includegraphics[scale=.8]{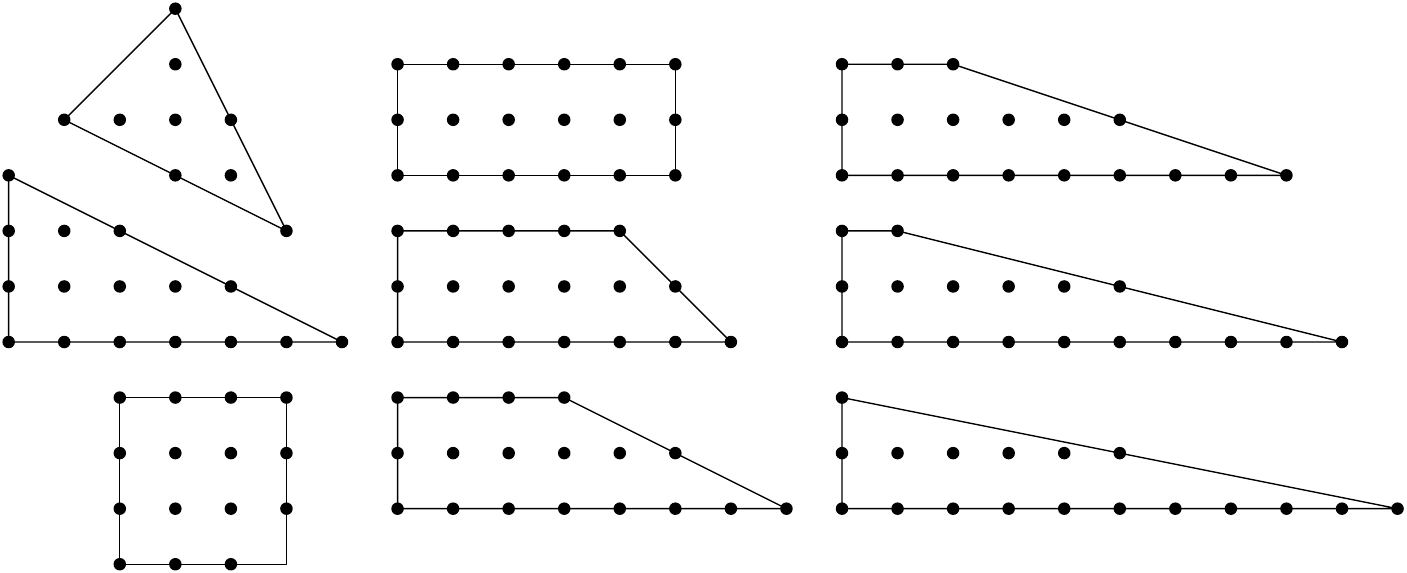}
\caption{The maximal polygons with $4$ interior lattice points.}
\label{fig:genus4_polygons}       
\end{figure}

Even when we have restricted to maximal polygons, there are two different flavors of polygons:  the {\emph{hyperelliptic} polygons}\index{hyperelliptic polygon}, for which $P_{int}$ as a line segment, and the {\emph{nonhyperelliptic} polygons}\index{nonhyperelliptic polygon}, for which $P_{int}$ is a two-dimensional polygon.  See Figure \ref{fig:genus4_polygons} for all the maximal polygons with $4$ interior lattice points, up to equivalence.  The leftmost three are nonhyperelliptic, and the other six are hyperelliptic.

How do we know there aren't any other maximal polygons with $4$ interior lattice points?  For the hyperelliptic case, \cite{koelman}  classifies all maximal hyperelliptic curves: they are a family of trapezoids interpolating between a hyperelliptic rectangle and a hyperelliptic triangle (this result is also presented in \cite{movingout}).  For the nonhyperelliptic polygons, we have the following result.

\begin{proposition}[Lemma 2.2.13 in \cite{koelman}; also Theorem 5 in \cite{movingout}]\label{prop:pushout}  Let $P$ be a maximal nonhyperelliptic polygon, with $P_{int}$ its interior polygon.  Then $P$ is obtained from $P_{int}$ by ``pushing out'' the edges of $P_{int}$.  More formally, if $P_{int}=\bigcap_{i=s}^s H_i$, where $H_i$ is the half-plane defined by the inequality $a_ix+b_iy\leq c_i$ (with $a_i,b_i,c_i$ relatively prime integers), then $P_{int}=\bigcap_{i=1}^s H'_i$, where $H'_i$ is the half-plane defined by the inequality $a_ix+b_iy\leq c_i+1$.
\end{proposition}

To find all nonhyperelliptic lattice polygons with $g$ interior lattice points, it thus suffices to find all lattice polygons with $g$ lattice points total, and then to determine which can be pushed out to form a lattice polygon.

\begin{exercise}  Using Proposition \ref{prop:pushout}, verify that Figure \ref{fig:genus4_polygons} does indeed contain all maximal nonhyperelliptic polygons with $4$ interior lattice points.  Then find all maximal nonhyperelliptic polygons with $5$ interior lattice points.
\end{exercise}

\begin{exercise}\label{exercise:hyperelliptic}  Determine which troplanar graphs of genus $g$ come from hyperelliptic Newton polygons.  (\textbf{Hint}:  if $g=3$, there are three such graphs, namely the middle three graphs from Figure \ref{fig:genus3_graphs}.)
\end{exercise}

\resproject{Study the properties of lattice polygons, stratified by the number of interior lattice points $g$.  (A great starting point for exploring these topics are the papers \cite{movingout} and \cite{cv}.)  For example:
 Given a maximal polygon $P$, let $n(P)$ be the number of subpolygons of $P$ with the same set of interior lattice points.   For which polygons is $n(P)$ equal to $1$?  What upper bounds can we find on $n(P)$, in terms of $g$?  How big is $n(P)$ on average?  (This gives us an idea of how much time we save by considering only maximal polygons when studying troplanar graphs.)
}

\begin{example} \label{example:genus3} Let us find all troplanar graphs of genus $3$.  (This will mirror arguments found in \cite{blmpr} and \cite{bjms}.) There are exactly five trivalent connected graphs of genus $3$ \cite{balaban}, namely those appearing in Figure \ref{fig:genus3_graphs}.  By Exercise \ref{exercise:skeleton_properties}, these are the only possible graphs that could be troplanar.  We now must which determine which of the five are actually achievable.

\begin{figure}[hbt]
\sidecaption
\includegraphics[scale=.65]{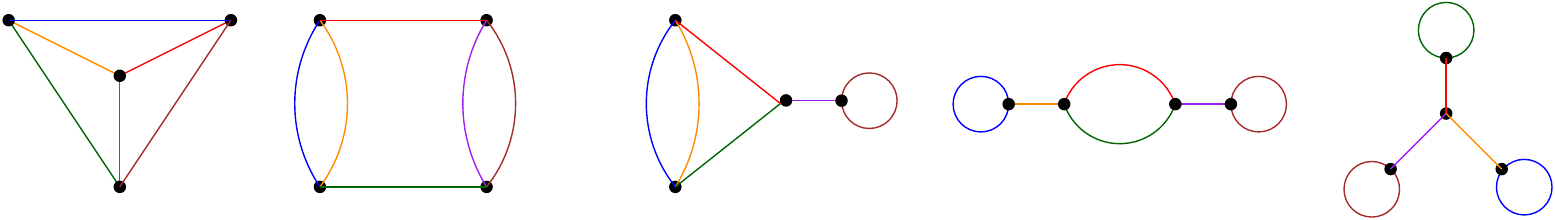}
\caption{The five candidate graphs of genus $3$.}
\label{fig:genus3_graphs}       
\end{figure}

Let us determine which Newton polygons are possible.  As mentioned previously, it suffices to take $P$ maximal.  We will focus on nonhyperelliptic polygons; the hyperelliptic ones are covered by Exercise \ref{exercise:hyperelliptic}.  It turns out that the only nonhyperelliptic polygon with $3$ interior lattice points, up to equivalence, is $T_4$, the triangle of degree $4$.  This is because the only lattice polygon (again, up to equivalence) with three lattice points is the triangle of degree $1$, which pushes out to $T_4$.   Figure \ref{fig:genus3_triangulations} shows triangulations of $T_4$ that give tropical curves whose skeletons are the first four graphs from Figure \ref{fig:genus3_graphs}, so we know that those four graphs are all troplanar.

\begin{figure}[hbt]
\sidecaption
\includegraphics[scale=.65]{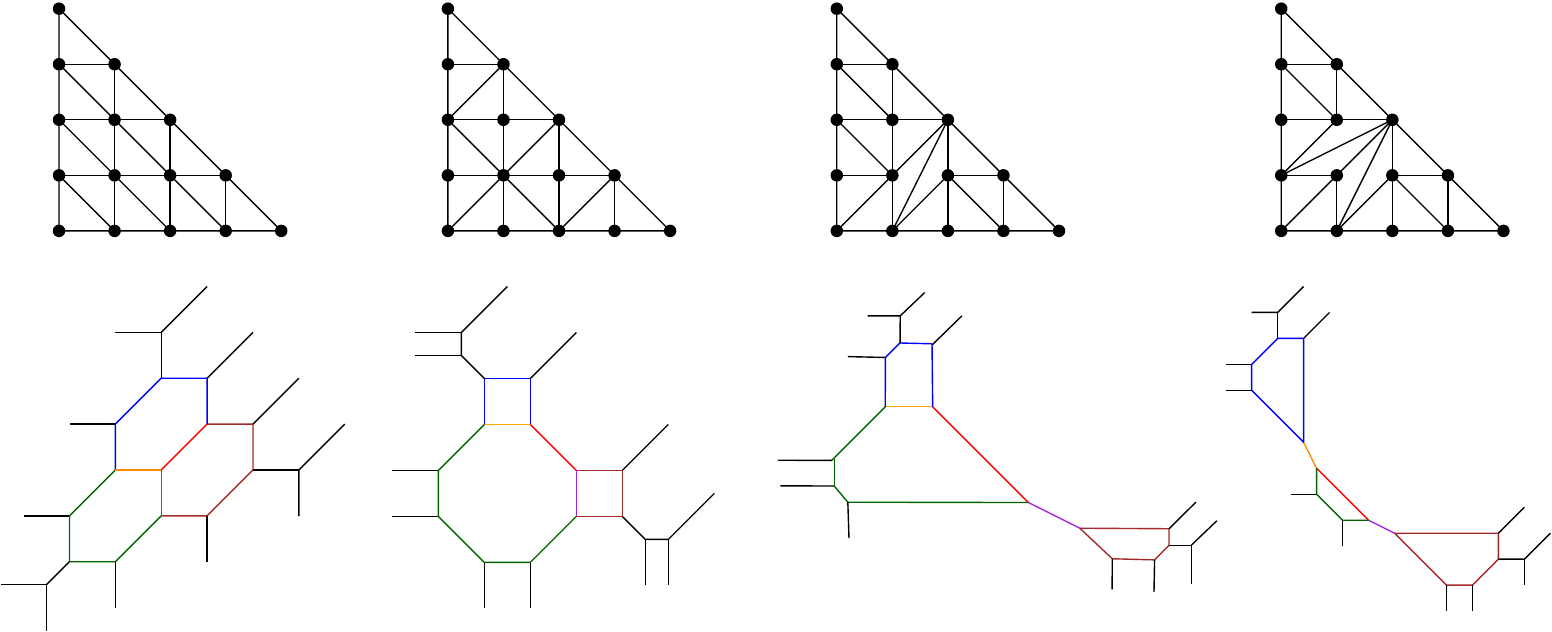}
\caption{Four triangulations, giving us four tropical curves whose skeletons are the first four graphs in Figure \ref{fig:genus3_graphs}.}
\label{fig:genus3_triangulations}       
\end{figure}

Lets now argue that the fifth graph, sometimes called the lollipop graph of genus $3$, is not troplanar.  Note that any bridge\footnote{A \emph{bridge} in a connected graph is an edge that, if removed from the graph, would disconnect the graph.} in troplanar graph must be dual to a \emph{split} in the subdivision of $T_4$, which is an edge goes from one boundary point to another, with some interior lattice points on each side and none in the edge's interior.  So, any triangulation of the triangle of degree $4$ that gives us the lollipop graph would have three splits.  All possible splits in the triangle are illustrated in Figure \ref{fig:splits}; however, no more than two of them can coexist in the same triangulation due to intersections, meaning we cannot obtain the lollipop graph.  We conclude that there are four troplanar graphs of genus $3$:  the first four graphs in Figure \ref{fig:genus3_graphs}.

\begin{figure}[hbt]
\sidecaption
\includegraphics[scale=.65]{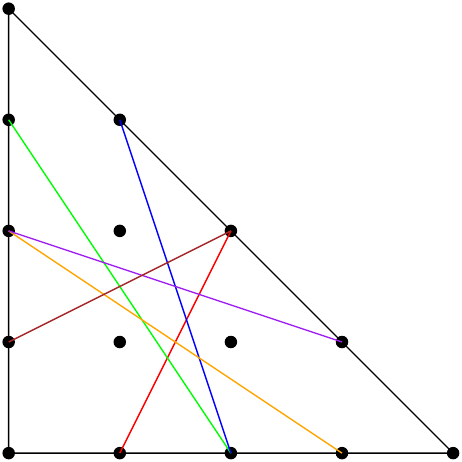}
\caption{Twelve splits, any three of which have at least one intersection point away from the boundary.}
\label{fig:splits}       
\end{figure}

\end{example}

\begin{figure}[hbt]
\sidecaption
\includegraphics[scale=.65]{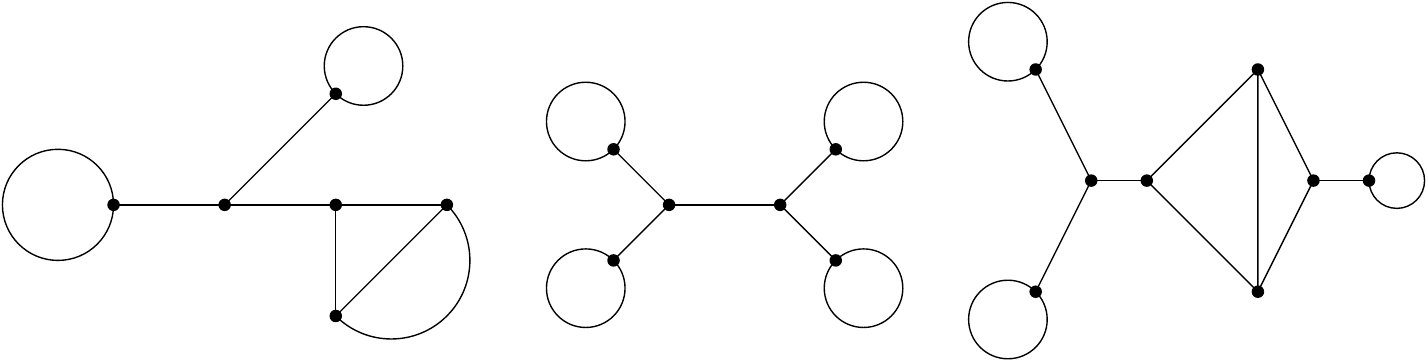}
\caption{Three sprawling graphs.  Note that the vertex that disconnects the graph into three pieces need not be unique.}
\label{fig:sprawling}       
\end{figure}

The fact that the lollipop graph did not appear also follows from a more general result about structures that cannot appear in troplanar graphs.  We say a connected, trivalent graph is {\emph{sprawling}}\index{sprawling graph} if removing a single vertex splits the graph into three pieces.  Several examples of sprawling graphs appear in Figure \ref{fig:sprawling}.  

\begin{proposition}[Proposition 4.1 in \cite{sprawling}]  \label{prop:sprawling} A sprawling graph cannot be troplanar.
\end{proposition}

Although this result was originally proved in \cite{sprawling}, the ``sprawling'' terminology comes from \cite{blmpr}, which offers an alternate proof.

\chproblem{  Prove Proposition \ref{prop:sprawling}.  (\textbf{Hint:}  Consider the structure of the dual triangulation of a smooth tropical curve with a sprawling skeleton.)
}

\begin{figure}[hbt]
\sidecaption
\includegraphics[scale=.5]{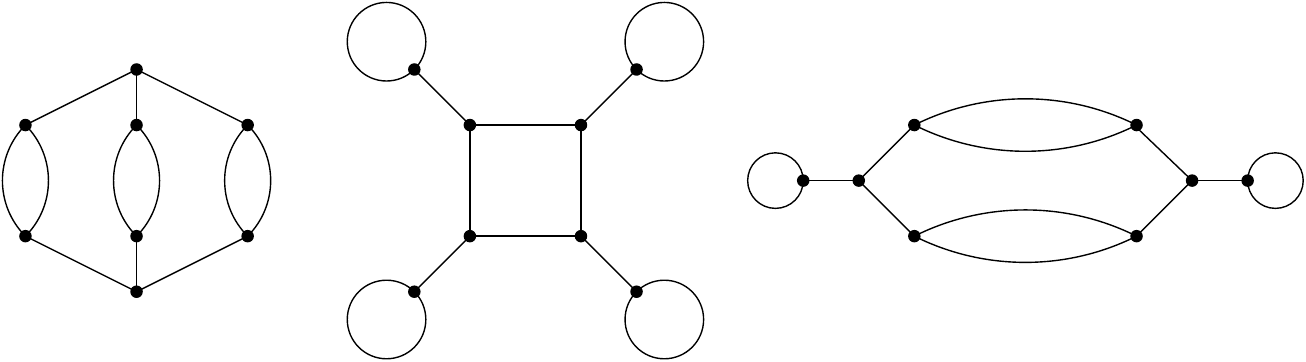}
\caption{Three graphs of genus $5$ that are \emph{not} troplanar}
\label{fig:genus5}       
\end{figure}

\chproblem{Show that the graphs in Figure \ref{fig:genus5} are not troplanar.}

\resproject{Find ``forbidden structures'' that never appear in troplanar graphs.  (Proposition \ref{prop:sprawling} gives an example of such a forbidden structure.  Another is given in \cite{morrison-hyperelliptic}.)}

\chproblem{There are $17$ trivalent connected graphs of genus $4$ \cite{balaban}.  Determine which of them are troplanar.  Note that the only Newton polygons you need to consider are those illustrated in Figure \ref{fig:genus4_polygons}.  (If you've already done Exercise \ref{exercise:hyperelliptic}, you can ignore six of the polygons!)}

In general, counting the number of tropically planar graphs of genus $g$ can be accomplished as follows:
\begin{enumerate}
\item  Find all maximal lattice polygons $P$ with $g$ interior lattice points, perhaps following \cite{movingout}.
\item  Find all regular unimodular triangulations of each $P$ from step 1, perhaps with \texttt{polymake} or \texttt{TOPCOM}.
\item  Find the dual skeletons to the triangulations from step 2, and sort them into isomorphism classes.
\end{enumerate}
This algorithm was implemented in \cite{bjms}, and was used to determine that the numbers of troplanar graphs of genus $2,3,4,$ and $5$ are $2$, $4$, $13$, and $37$, respectively.  This was pushed further as part of the Williams SMALL 2017 REU to genus $6$ ($151$ troplanar graphs) and genus $7$ ($672$ troplanar graphs).

\resproject{Find a more efficient way to determine the number of troplanar graphs of genus $g$ than the algorithm outlined above.}

\resproject{Study how the number of troplanar graphs of genus $g$ grows with $g$.  Can you find upper and lower bounds?  Can you determined its asymptotic behavior?  (Preliminary work in this direction was done in the Williams College SMALL REU in 2017.)}

So far we have considered skeletons from a purely combinatorial perspective.  Now we include the data of lengths on each edge of the graph, giving us a {\emph{metric graph}}\index{metric graph}.  A natural impulse is to sum up all the Euclidean lengths of the edges of the embedded tropical curve that make up a given edge of the skeleton, and declare that to be its length.  Unfortunately this definition of length is not invariant under the natural transformations that we apply to our Newton polygons.  This leads us to use the following definition. 

\begin{definition}  Let $P_1,P_2\in\mathbb{R}^2$ be distinct points such that the line segment $\overline{P_1P_2}$ has rational slope (or is vertical).  Write the vector from $P_1$ to $P_2$ as $\lambda\times \left<a,b\right>$, where $a,b\in\mathbb{Z}$ with $\gcd(a,b)=1$ and $\lambda\in\mathbb{R}^+$.  The {\emph{lattice length}}\index{lattice length} of the line segment $\overline{P_1P_2}$ is defined to be $\lambda$.
\end{definition}

When considering a tropical plane curve, we measure the lengths of its finite edges by lattice length.  These lengths are then added up appropriately to assign lengths to the edges of the skeleton.

\begin{example}\label{example:metric}  Consider the tropical plane curve illustrated on the top in Figure \ref{fig:lattice_length_example}.  Below it is the collection of all bounded edges in the curve, labelled with their lattice lengths.  As pictured, the skeleton is a graph consisting of two vertices joined by an edge, with a loop attached to each vertex.  The length of the middle edge in the skeleton is $1$; the lengths of the loops are $2+1+1+3+5=12$ and $6+3+3+1+1+1=15$.  (Note that one bounded edge from the tropical curve does not contribute to the skeleton.)

\begin{figure}[hbt]
\sidecaption
\includegraphics[scale=.65]{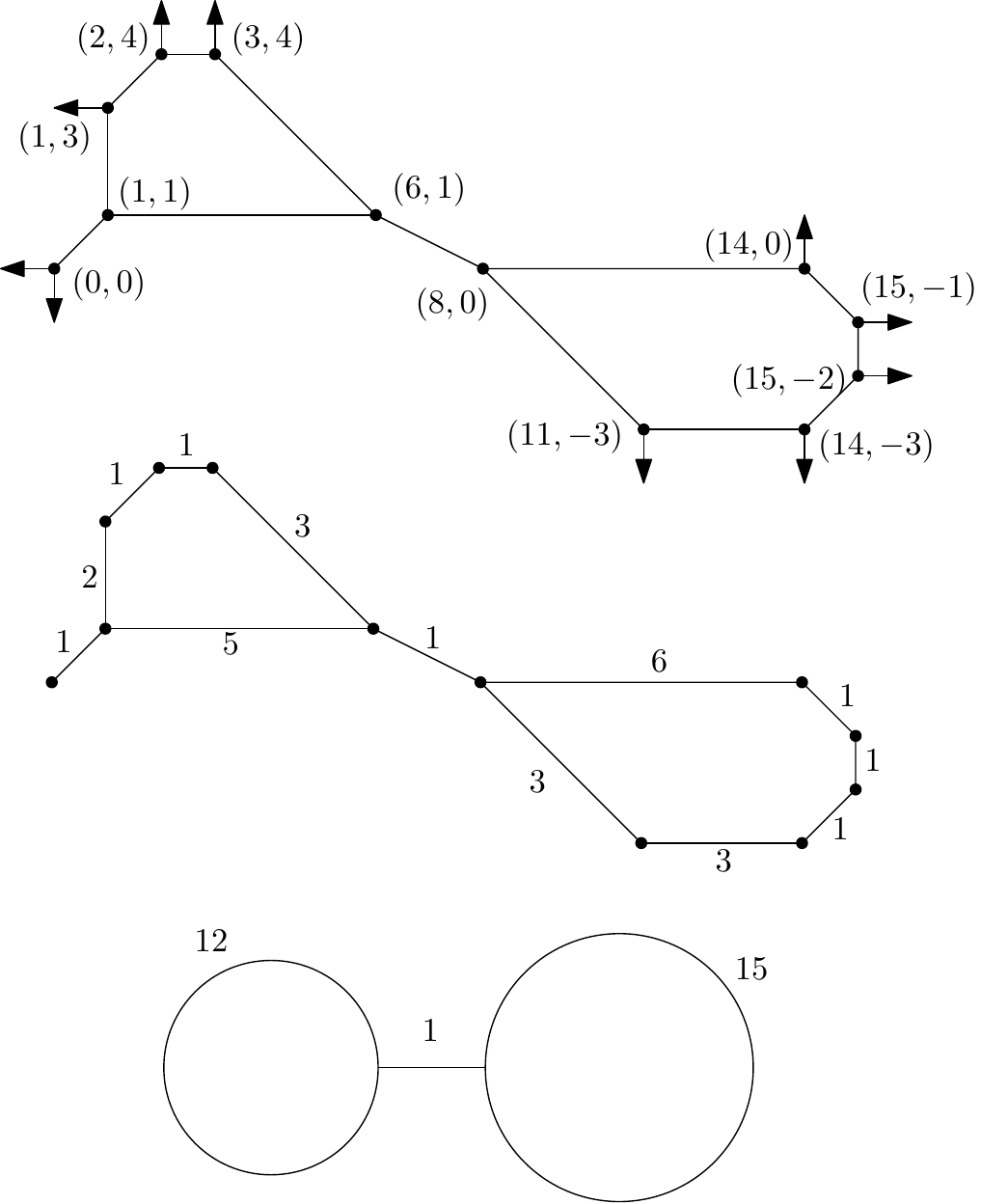}
\caption{A tropical curve with lattice lengths labelled, and the resulting lengths on the skeleton}
\label{fig:lattice_length_example}       
\end{figure}

\end{example}

When we say that a metric graph is troplanar, we mean that it is the skeleton of a smooth tropical plane curve \emph{giving those edge lengths}.  So the metric graph at the bottom of Figure \ref{fig:lattice_length_example} is troplanar.

\chproblem{Let $P$ be a $2\times 3$ lattice rectangle.  Find all troplanar metric graphs that are the skeleton of a smooth tropical curve with that Newton polygon.  (\textbf{Hint}: in some sense you can get most, but not all, graphs of genus $2$.)}

The algorithm presented in \cite{blmpr} did not simply find the combinatorial types of troplanar graphs; it computed, up to closure, all \emph{metric} graphs of genus at most $5$ that appeared as the skeleton of a smooth tropical plane curve.  In their Theorem 5.1, they use this computation to characterize exactly which metric graphs of genus $3$ are troplanar.  Beyond the lollipop graph not appearing (regardless of the edge lengths), there are nontrivial edge length restrictions on the other four combinatorial types of graphs.  Rather than presenting their full result here, we give a consequence of it.

\begin{theorem}[Corollary 5.2 in \cite{bjms}]  Approximately $29.5\%$ of all metric graphs of genus $3$ are troplanar.
\end{theorem}
 This probability is computed by considering the \emph{moduli space of graphs of genus $3$} \cite{bmv,chan}.  This is a six-dimensional space, corresponding to the six edges a trivalent graph of genus $3$ has.  This space is not compact, since edge lengths can be arbitrarily long; so consider the subspace consisting of graphs with total length equal to $1$; up to scaling, every metric graph can be represented in this way. Give each of the five combinatorial types of graphs (as illustrated in Figure \ref{fig:genus3_graphs}) an equal weight, and compute the volume of the space of troplanar graphs within this $5$-dimensional space.  This computation gives about $0.295$, or $29.5\%$.  
 
 \chproblem{Show that neither of the metric graphs illustrated in Figure \ref{fig:nontroplanar} are troplanar.  (This follows from the characterization given in \cite[Theorem 5.1]{bjms}; try to give your own argument.)
 
 \begin{figure}[hbt]
\sidecaption
\includegraphics[scale=.60]{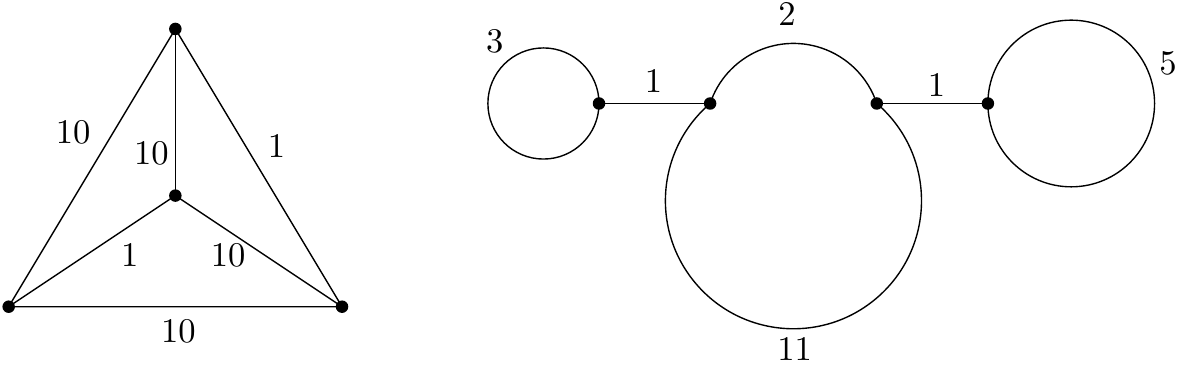}
\caption{Two metric graphs that aren't troplanar}
\label{fig:nontroplanar}       
\end{figure}

 }

\resproject{Determine which metric graphs arise as the skeleton of a smooth tropical plane curve, perhaps under certain restrictions.  For instance:
\begin{itemize}
\item  Characterize exactly which metric graphs arise from hyperelliptic polygons, as explored in \cite{morrison-hyperelliptic}.
\item  Characterize which metric graphs arise from \emph{honeycomb polygons}, a key tool in \cite{bjms}.
\item  Characterize which metric graphs are troplanar with as many degrees of freedom as possible on their edge lengths.  In \cite{bjms}, this maximum number of degrees of freedom was shown to be $2g+1$, at least for $g\geq 8$.
\end{itemize}}

All of our questions have been posed for smooth tropical plane curves.  Of course, we can also consider tropical curves with singularities.  We say a tropical curve is \emph{nodal} if, in the dual subdivision, all polygons besides the triangles of area $1/2$ are quadrilaterals of area $1$.  A vertex in a nodal tropical curve dual to such a quadrilateral is called a \emph{node}.

\begin{example}\label{example:nodal}  Figure \ref{fig:nodal} presents an example of a nodal tropical curve with its dual Newton subdivision.  We can still consider a skeleton of the curve by interpreting each nodal crossing in the tropical curve as two edges in the graph that happen to look like they're crossing.  The resulting skeleton is pictured on the right.

 \begin{figure}[hbt]
\sidecaption
\includegraphics[scale=1]{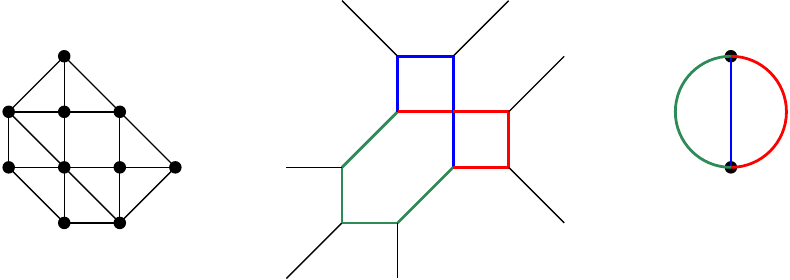}
\caption{A nodal tropical curve and its skeleton}
\label{fig:nodal}       
\end{figure}
\end{example}

It was shown in \cite{sprawling} that \emph{every} connected trivalent graph can be realized in a nodal tropical plane curve.  Given a connected trivalent graph $G$, let $N(G)$ be the {\emph{tropical crossing number}}\index{tropical crossing number} of $G$, which is the smallest number of nodes required to achieve $G$ as the skeleton of a nodal tropical curve.  For instance, $N(G)=0$ if and only if $G$ is troplanar.

\resproject{
Study the tropical crossing number.  Can you determine its value explicitly for certain families of graphs?  (Note that if this question is being posed for metric graphs, $N(G)$ does depend on the edge lengths.)
}

\section{Tropical Geometry in Three Dimensions}
\label{section:three_dimensions}

Moving beyond the plane into three-dimensional space, we consider tropical polynomials in three variables $x,y,$ and $z$. Such a polynomial can be written as
\begin{equation}p(x,y,z)=\bigoplus_{(i,j,k)\in S} c_{ijk}\odot x^i\odot y^j\odot z^k,\end{equation}
where $S$ is the set of all exponent vectors that appear in $p(x,y,z)$.  This polynomial defines a {\emph{tropical surface}}\index{tropical surface}, the set of all points in $\mathbb{R}^3$ where the maximum defined by the polynomial is achieved at least twice. Again, we denote this tropical surface $\mathcal{T}(p)$.

\begin{figure}[hbt]
\sidecaption
\includegraphics[scale=.65]{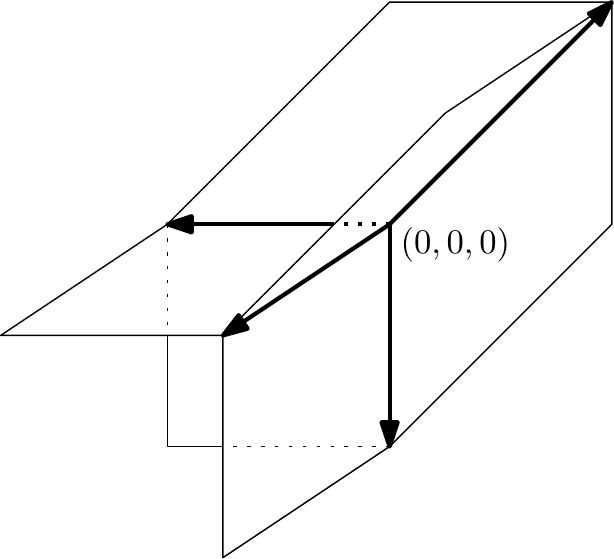}
\caption{The tropical plane defined by $x\oplus y\oplus z\oplus 0$}
\label{fig:plane}       
\end{figure}

\begin{example}\label{example:plane} Let $p(x,y,z)=x\oplus y\oplus z\oplus 0$. The tropical surface $\mathcal{T}(p)$ is illustrated in Figure  \ref{fig:plane}.  It consists of the origin $(0,0,0)$; four rays, pointing in the directions $\left<-1,0,0\right>$, $\left<0,-1,0\right>$, $\left<0,0,-1\right>$, and $\left<1,1,1\right>$; and six two-dimensional pieces, each obtained as the positive linear span of two of the rays.  Such two-dimensional pieces of a tropical surface are called \emph{two-dimensional cells}.  Because of the form of $p(x,y,z)$, we call $\mathcal{T}(p)$ a {\emph{tropical plane}}\index{tropical plane}.
\end{example}

\subsection{Tropical Surfaces and the Duality Theorem}

The Duality Theorem still holds for tropical polynomials in three variables and the surfaces they define\footnote{Indeed, a Duality Theorem holds for all tropical varieties defined by a single equation in any number of variables; see {\cite[Proposition 3.1.6]{maclagan-sturmfels}}.}.  This time, instead of a Newton polygon we consider a \emph{Newton polytope}, the convex hull of all exponent vectors appearing in the polynomial. (We will assume that the Newton polytope is three-dimensional to avoid certain degenerate cases.)  To find an induced subdivision, we again associate heights to each lattice point of the Newton polytope; this time, however, we must compute our upper convex hull in four-dimensional space.  We then have the following correspondence between parts of the tropical surface $S=\mathcal{T}(p)$ and the subdivision of $\textrm{Newt}(p)$:
\begin{itemize}
\item  Vertices in $S$ correspond to $3$-dimensional polytopes in the subdivision.
\item  Rays in $S$ correspond to boundary two-dimensional faces.
\item  Edges in $S$ correspond to interior two-dimensional faces.
\item  Unbounded two-dimensional cells in $S$ correspond to boundary edges.
\item  Bounded two-dimensional cells in $S$ correspond to interior edges.
\end{itemize}
As was the case for tropical plane curves, the relationships and geometry of all these pieces of the tropical surface are dictated by the subdivision.  For instance, two vertices are joined by an edge if and only if the corresponding polytopes share a face; and that edge is perpendicular to the shared face.

We say that a subdivision of a polytope is a \emph{unimodular tetrahedralization} if all polytopes in the subdivision are tetrahedra of volume $\frac{1}{6}$, which is the smallest possible volume.  We say that a tropical surface $\mathcal{T}(p)$ is \emph{smooth} if the induced subdivision of $\textrm{Newt}(p)$ is a unimodular tetrahedralization.   If $\textrm{Newt}(p)$ is the tetrahedron with vertices at $(0,0,0)$, $(d,0,0)$, $(0,d,0)$, and $(0,0,d)$, we say that $p(x,y,z)$ \emph{has degree $d$}.

\example{\label{example:almost_cube}Let 
\begin{equation}
f(x,y,z)=(x y\odot z)\oplus (-42\odot x y)\oplus x\oplus y\oplus z\oplus (-42),\end{equation} and let $P=\text{Newt}(f)$.  The polytope $P$ looks like a cube with two tetrahedra sliced off, as illustrated to the left in Figure \ref{fig:subdivided_polytope}.  Every term has coefficient $0$, except for the $(0,0,0)$ and $(1,1,0)$ terms, which have a very negative coefficient.  This means that in the subdivision, we will end up with two smaller tetrahedra with vertices at $(1,0,0)$, $(0,1,0)$, $(1,1,0)$, and $(1,1,1)$; and at $(0,0,0)$, $(0,0,1)$, $(1,0,0)$, and $(0,1,0)$; as well as a larger tetrahedron at $(0,0,1)$, $(1,0,0)$, $(0,1,0)$, and $(1,1,1)$\footnote{To prove this rigorously, we would need to show that the hyperplane in $\mathbb{R}^4$ containing the points $(1,0,0,0)$, $(0,1,0,0)$, $(1,1,0,0)$, and $(1,1,1,-42)$ lies strictly above the points $(0,0,1,0)$ and $(0,0,0,-42)$; as well as two other similar such statements, one for each of the other tetrahedra.  (In fact, the hyperplane we get from the middle tetrahedron in $(x,y,z,w)$-space is just defined by $w=0$, and certainly the other two lifted points $(0,0,0,-42)$ and $(1,1,0,-42)$ lie below this hyperplane.)}.  This is illustrated in Figure \ref{fig:subdivided_polytope}.

\begin{figure}[hbt]
\sidecaption
\includegraphics[scale=1]{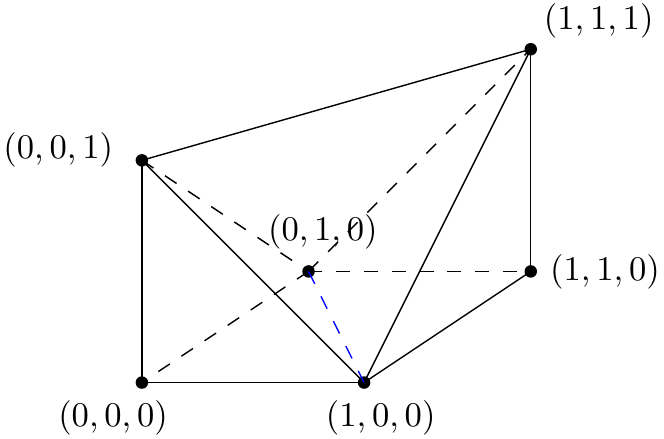}
\caption{The subdivided Newton polytope from Example~\ref{example:almost_cube}}
\label{fig:subdivided_polytope}       
\end{figure}

The tropical surface $\mathcal{T}(f)$ has three vertices, corresponding to the three tetrahedra.  We can find their coordinates by computing the four-way ties.
 \begin{itemize}
 \item  From $-42=x=y=z$, we have a vertex at $(-42,-42,-42)$.
 \item  From $x=y=z=x+y+z$, we have a vertex at $(0,0,0)$.
 \item  From $x=y=-42+x+y=x+y+z$, we have a vertex at $(42,42,-42)$.
 \end{itemize}
 The vertex at $(0,0,0)$ connects to the other two vertices by a line segment.  The vertex $(-42,-42,-42)$ will have three rays, pointing in the directions $\left<-1,0,0\right>$,$\left<0,-1,0\right>$, and $\left<0,0,-1\right>$.  The vertex $(0,0,0)$ will have two rays, pointing in the directions $\left<1,-1,1,\right>$ and $\left<-1,1,1,\right>$.  Finally, the vertex $(42,42,-42)$ will  have three rays, pointing in the directions  $\left<1,0,0\right>$,$\left<0,1,0\right>$, and $\left<0,0,1\right>$.  Ignoring the two-dimensional pieces, our tropical surface looks as pictured in Figure \ref{fig:tropical_surface}.

\begin{figure}[hbt]
\sidecaption
\includegraphics[scale=1]{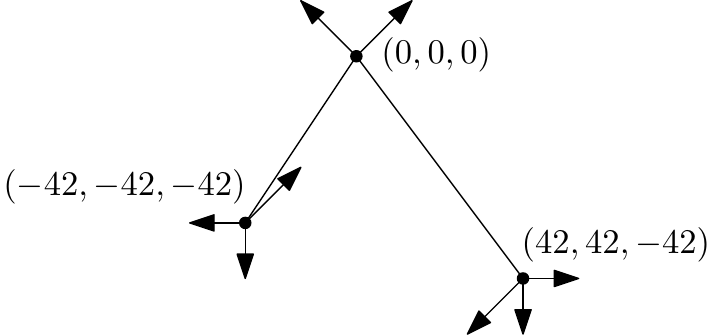}
\caption{The one-dimensional pieces of the surface from from Example~\ref{example:almost_cube}}
\label{fig:tropical_surface}       
\end{figure}

 We fill in two-dimensional pieces between adjacent rays and edges.  This will give a total of $12$ unbounded $2$-dimensional pieces, corresponding to the $12$ edges in our tetrahedralization.  All are unbounded, since all edges in the tetrahedralization are exterior.
 
Note that this tropical surface is \emph{not} smooth.  Even though our induced subdivision is a tetrahedralization, it is not unimodular since the tetrahedra don't all have volume $\frac{1}{6}$: the tetrahedron in the middle has volume $1/3$.
 
.}

\chproblem{\label{cp:surfaces} Show that the tropical polynomial of degree $2$ defined by

\begin{equation}
\begin{matrix}
f & = &  (-3\odot{x}^{2})\oplus(-4\odot xy)\oplus xz\oplus (-7\odot{y}^{2})\oplus(-2\odot yz)\\ & &\oplus (-1\odot{z}^{2})\oplus
x\oplus y\oplus (-2\odot z)\oplus (-7)
\end{matrix}
\end{equation}
is a smooth tropical surface.  Determine how many vertices, edges, rays, bounded two-dimensional cells, and unbounded two-dimensional cells there are.  Do the same for the tropical polynomial of degree $3$ defined by
\begin{equation}
\begin{matrix}
g & = &
(-23\odot {x}^{3})\oplus(-15\odot{x}^{2}y)\oplus(-7\odot{x}^{2}z)\oplus(-15\odot x{y}^{2})\oplus xyz\oplus (-3\odot x{
z}^{2})
\\ & &\oplus (-25\odot {y}^{3})\oplus (-6\odot {y}^{2}z)\oplus (-10\odot y{z}^{2})\oplus (-20\odot {z}^{3})\oplus   (-2\odot 
{x}^{2})
\\ & &\oplus (-6\odot xy)\oplus (-1\odot xz)\oplus (-14\odot {y}^{2})\oplus yz\oplus (-9\odot {z}^{2})\oplus (-11\odot x)
\\ & &\oplus 
(-4\odot y)\oplus (-9\odot z)\oplus (-21).
\end{matrix}
\end{equation}
You will almost certainly want to use a computer to help with this!  After you try this Challenge Problem, you should check your counts against the following theorem.  
}

\begin{theorem}[Theorem 4.5.2 in \cite{maclagan-sturmfels}]  A smooth tropical surface of degree $d$ has
\begin{itemize}
\item $d^3$ vertices,
\item  $2d^2(d-1)$ edges,
\item  $4d^2$ rays,
\item $d(d-1)(7d-11)/6$ bounded two-dimensional cells,
\item  $6d^2$ unbounded two-dimensional cells.
\end{itemize}
Its Euler characteristic\footnote{Intuitively, this is the number of bounded regions of $\mathbb{R}^3$ encapsulated by part of the surface.} is $\frac{(d-1)(d-2)(d-3)}{6}+1$.
\end{theorem}

\resproject{Study the geometry of smooth tropical surfaces.  For instance:

\begin{itemize}
\item A smooth surface of degree $3$ has $10$ bounded two-dimensional cells, each of which is a polygon, say with $n_i$ sides for the $i^{th}$ polygon.  What are the possible values for $n_1,\ldots,n_{10}$?  How can these $10$ polygons be arranged relative to each other?
\item A smooth surface of degree $4$ has Euler characteristic $1$, and so contains one polytope bounding a three-dimensional region.  Can we characterize which polytopes are possible?  (How many faces, how many edges, etc.)
\item  Moving on to smooth surfaces of degree grater than $4$, which have Euler characteristic greater than $1$, there are multiple polytopes that are part of the surface.  How can these polytopes be arranged?  (This is the surface analog of asking what the skeleton of a smooth tropical plane curve can be.)
\end{itemize}}

\subsection{Tropical curves in $\mathbb{R}^3$}

In usual geometry, if we intersect a pair of two-dimensional surfaces in $\mathbb{R}^3$, we expect to get a one-dimensional curve.  This also holds in tropical geometry, if we are willing to assume stable intersections to avoid the overlap of two-dimensional pieces.  There is still a duality theorem for tropical curves in $\mathbb{R}^3$ that arise as the intersection of two tropical surfaces, although it requires a bit more machinery.

Given two lattice polytopes $P,Q\subset\mathbb{R}^n$, place $P$ and $Q$ in $(n+1)$-dimensional space by giving every point in $P$ an extra coordinate of $0$ and every coordinate of $Q$ an extra coordinate of $1$.   The {\emph{Cayley polytope}}\index{Cayley polytope} of $P$ and $Q$, written $\textrm{Cay}(P,Q)$, is the convex hull in $\mathbb{R}^{n+1}$ of this arrangement.

\begin{example}\label{example:cayley}  If $P=Q=\textrm{conv}(\{(0,0,0),(1,0,0), (0,1,0),(0,0,1)\})$, then $\textrm{Cay}(P,Q)$ is the convex hull of the eight points $(0,0,0,0)$, $(1,0,0,0)$, $(0,1,0,0)$, $(0,0,1,0)$, $(0,0,0,1)$, $(1,0,0,1)$, $(0,1,0,1)$, and $(0,0,1,1)$ in $\mathbb{R}^4$.
\end{example}

Suppose $p(x,y,z)$ and $q(x,y,z)$ are tropical polynomials in three variables, defining tropical surfaces $S_1$ and $S_2$, with intersection curve $C=S_1\cap_{st}S_2$.  Let $P=\textrm{Newt}(p)$ and $Q=\textrm{Newt}(q)$.  As we did with Newton polygons and Newton polytopes of single polynomials, we can find an induced subdivision of $\textrm{Cay}(P,Q)$. Each lattice point of $\textrm{Cay}(P,Q)$ is either a  lattice point of $P$ with an extra coordinate of $0$ or a lattice point of $Q$ with an extra coordinate of $1$, so we assign to each such lattice point a ``height'' based on the corresponding coefficient from the relevant polynomial.  We can then compute the induced subdivision of $\textrm{Cay}(P,Q)$ by looking at the upper convex hull in $\mathbb{R}^5$ of these lifted points.  This subdivision then splits $\textrm{Cay}(P,Q)$ into $4$-dimensional polytopes. Some of these polytopes have one vertex from $P$ and all others from $Q$, or vice versa; the other polytopes, with at least two vertices coming from each of $P$ and $Q$, are called the \emph{mixed cells} of the subdivision.  The duality theorem for complete intersection curves, stated fully in \cite[\S 4.6]{maclagan-sturmfels}, then says that the vertices of $C$ correspond to the mixed cells of this subdivision.

If all cells in the Cayley subdivision have the minimum possible volume (which turns out to be $1/24$), we call the tropical curve \emph{smooth}.  In this case it turns out that $P\cap_{st}Q=P\cap Q$.  We can still talk about the skeletons of tropical curves in $\mathbb{R}^3$, retracting rays and leaves to obtain the desired graph.  Again we still refer to the genus of the graph, although since it might not be a planar graph we need to define genus as $E-V+1$.

\begin{theorem}[Theorem 4.6.20 in \cite{maclagan-sturmfels}]\label{theorem:tropical_curve_parts}  Let $f(x,y,z)$ and $g(x,y,z)$ be tropical polynomials with degrees $d$ and $e$, respectively, such that $C=\mathcal{T}(f)\cap \mathcal{T}(g)$ is a smooth tropical curve.  Then $C$ has
\begin{itemize}
\item $d^2e+de^2$ vertices,
\item  $(3/2)d^2e+(3/2)de^2-2de$ edges,
\item  $4de$ rays, and
\item  genus equal to $(1/2)d^2e+(1/2)de^2-2de+1$.
\end{itemize}
\end{theorem}

\begin{example}\label{example:3d_line}
Let $p(x,y,z)=(-1\odot x)\oplus (-1\odot y)\oplus z\oplus 1$ and $q(x,y,z)=(-2\odot x)\oplus (1\odot y)\oplus (1\odot z)\oplus (-1)$.  Then $\textrm{Newt}(p)$ and $\textrm{Newt}(q)$ are $P$ and $Q$ from Example \ref{example:cayley}.  Using the \texttt{Macaulay2} package \texttt{Polyhedra}\footnote{The Polyhedra package defaults to the min convention rather than the max.  This means we have to negate all the coefficients before we find the decomposition.}, we compute the subdivision of $\textrm{Cay}(P,Q)$.  It consists of four cells:
\begin{equation}
\Delta_1=\textrm{conv}(\{(0, 0, 0, 0), (1, 0, 0, 0), (0, 1, 0, 0), (0, 0, 1, 0), (0, 1,0, 1)\})\end{equation}
\begin{equation}
\Delta_2=\textrm{conv}(\{(0, 0, 0, 0), (1, 0, 0, 0), (0, 0, 1, 0), (0, 1, 0, 1),(0, 0, 1, 1)\})\end{equation}
\begin{equation}
\Delta_3=\textrm{conv}(\{(0, 0, 0, 0), (1, 0, 0, 0), (1, 0, 0, 1), (0, 1,0, 1), (0, 0, 1, 1)
\})\end{equation}
\begin{equation}
\Delta_4=\textrm{conv}(\{(0, 0, 0, 0), (0, 0, 0, 1), (1, 0, 0, 1),(0, 1, 0, 1), (0, 0, 1, 1)))
\})\end{equation}
Each cell has volume $1/24$, so the tropical intersection curve is smooth; as the intersection of two tropical planes, we call it a \emph{tropical line in $\mathbb{R}^3$}.  Of the four cells, only $\Delta_2$ and $\Delta_3$ are mixed cells.  This means the line $P\cap Q$ has two vertices.  The vertex $(a,b,c)$ coming from $\Delta_2$ arises from a three-way tie between the $(0,0,0)$, $(1,0,0)$, $(0,0,1)$ terms of $p$ and a two-way tie between the $(0,1,0)$ and $(0,0,1)$ terms of $q$. Written in conventional notation, we have
$1=-1+a=c$, so $a=2$ and $c=1$.  We also have $1+b=1+c$, so $b=c=1$.  Thus there is a vertex at $(a,b,c)=(2,1,1)$.  In the next exercise, you'll find the other vertex, as well as the rest of the line.
\end{example}

\begin{exercise}  Draw the tropical line from the previous example.  Be sure to check your answer against Theorem \ref{theorem:tropical_curve_parts} with $d=e=1$.
\end{exercise}

\chproblem{Show that the tropical surfaces from Challenge Problem \ref{cp:surfaces} intersect in a smooth tropical curve.  Show that the skeleton of the curve is the complete bipartite graph $K_{3,3}$.}

\resproject{Which graphs of genus $4$ arise in smooth tropical curves that are the intersection of a tropical surface of degree $2$ and a tropical surface of degree $3$?  For instance, are any of these graphs sprawling?

More generally:  which graphs of genus $(1/2)d^2e+(1/2)de^2-2de+1$ arise as the skeleton of a smooth tropical curve that is the intersection of a surface of degree $d$ with a surface of degree $e$?

(You can approach these questions considering the graphs either combinatorially, or as metric graphs.)}

\resproject{Let $Q_1$ and $Q_2$ be two smooth tropical surfaces of degree $2$.  Study the possibilities of the intersection $Q_1\cap_{st}Q_2$, possibly through a similar lens as \cite{fno}.  (If the intersection is a smooth curve, then it has genus $1$ by Theorem \ref{theorem:tropical_curve_parts}, and we understand its combinatorial properties very well.  What other intersections are possible?)}

One noteworthy difference between classical geometry and tropical geometry is that in tropical geometry, not all planes look the same.  In the previous section, we studied tropical curves as a subset of the usual plane $\mathbb{R}^2$.  But this plane is combinatorially different from, say, the tropical plane from Example \ref{example:plane}.  A natural question is then whether or not there are ``tropical plane curves'' besides those we studied in Section \ref{section:plane}; that is, whether certain tropical skeletons appear on tropical planes in $\mathbb{R}^3$ that did not arise from tropical curves in $\mathbb{R}^2$.  (We could ask the same for $2$-dimensional tropical planes in $\mathbb{R}^4$, or $\mathbb{R}^5$, or in general $\mathbb{R}^n$.)

Recent work shows that the answer is yes!  Recall that only $29.1\%$ of all graphs of genus $3$ appear in tropical curves in $\mathbb{R}^2$.   It is shown in \cite{hmrt} that \emph{every} metric graph of genus $3$, besides a family of measure zero, appears as a tropical curve in a tropical plane in $\mathbb{R}^3$, $\mathbb{R}^4$, or $\mathbb{R}^5$.  For example, they show that the lollipop graph appears as a tropical curve on a tropical plane in $\mathbb{R}^5$.  It is not known if their result is sharp; for instance, it is an open question if there are any graphs of genus $3$ that do not appear on a tropical plane in $\mathbb{R}^3$.

\resproject{Can the lollipop graph be realized on a tropical plane in $\mathbb{R}^3$ or $\mathbb{R}^4$?   More generally, which graphs can be realized on a tropical plane in $\mathbb{R}^n$, for different values of $n$?}

\section{Tropicalization}
\label{section:tropicalization}

 In this section we present the connections between \emph{algebraic geometry}, which studies solutions to usual polynomial equations, and tropical geometry, which studies solutions to tropical polynomial equations.  See \cite{maclagan-sturmfels} for a more complete treatment of this connection, and \cite{clo} for an undergraduate introduction to algebraic geometry.

Let $k$ be a field, and let $k[x_1,\ldots,x_n]$ be the polynomial ring in $n$ variables over $k$.  For an ideal $I\subset k[x_1,\ldots,x_n]$, the {\emph{affine variety}}\index{affine variety} defined by $I$ is
\begin{equation}
\textbf{V}(I)=\{\,(a_1,\ldots,a_n)\,|\,f(a_1,\ldots,a_n)=0\textrm{ for all $f\in I$}\,\}\subset k^n.\end{equation}
Given $f_1,\ldots,f_s\in k[x_1,\ldots,x_n]$, we can also define
\begin{equation}
\textbf{V}(f_1,\ldots,f_s)=\{\,(a_1,\ldots,a_n)\,|\,f_i(a_1,\ldots,a_n)=0\textrm{ for all $i$}\,\}\subset k^n.\end{equation}
If $I=\left<f_1,\cdots,f_s\right>$, then $\textbf{V}(I)=\textbf{V}(f_1,\ldots,f_s)$.  By Hilbert's Basis Theorem \cite{hilbert}\footnote{For a presentation in English, see \cite[\S2.5]{clo}.} every ideal in $k[x_1,\ldots x_n]$ has a finite set of generators,  so these two characterizations of affine varieties are equivalent.

Sometimes it is useful to work within the ambient space of the \emph{algebraic torus} $(k^*)^n$, where $k^*=k\setminus\{0\}$.  To do this we can let our ideal $I$ be a subset of $k[x_1^{\pm 1},\ldots,x_n^{\pm 1}]$, so that $\textbf{V}(I)\subset (k^*)^n$.

\subsection{Fields with valuation}

We will work with fields with an additional structure called a {\emph{valuation}}\index{valuation}.  A valuation on a field $k$ is a function  $\textrm{val}:k\rightarrow(\mathbb{R}\cup\{\infty\})$ such that
\begin{itemize}
\item $\textrm{val}(a)=\infty$ if and only if $a=0$.
\item $\textrm{val}(ab)=\textrm{val}(a)+\textrm{val}(b)$.
\item  $\textrm{val}(a+b)\geq\min\{\textrm{val}(a),\textrm{val}(b)\}$ with equality if $\textrm{val}(a)\neq \textrm{val}(b)$.
\end{itemize}
Every field has an example of a valuation called the \emph{trivial valuation}, defined by $\textrm{val}(0)=\infty$  and $\textrm{val}(a)=0$ for all $a\neq 0$.  Let's find some nontrivial valuations.

\begin{exercise}  Let $\mathbb{Q}$ be the field of rational numbers, and let $p$ be a prime number.  Define the {\emph{$p$-adic valuation}}\index{$p$-adic valuation} on $\mathbb{Q}$ by
\begin{equation}
\textrm{val}_p\left(p^k\frac{a}{b}\right)=k,
\end{equation}
where $a$ and $b$ are integers that aren't divisible by $p$.  Show that this is a valuation on~$\mathbb{Q}$.
\end{exercise}

\resproject{Study the sequences obtained by applying $p$-adic valuations to sequences of integers.  For instance, applying the $2$-adic valuation to the sequence of Fibonacci numbers
\begin{equation}
1,1,2,3,5,8,13,21,34,55,\dots
\end{equation}
gives the sequence
\begin{equation}
0,0,1,0,0,3,0,0,1,0,\dots
\end{equation}
We can think of this as \emph{tropicalizing} sequences of integers.  See \cite{amm, lengyel, fibonacci} for work done in this direction.
 }

\begin{exercise}  Let $K$ be a field and let $K((t))$ be the field of {\emph{Laurent series}}\index{Laurent series} over $K$, the nonzero elements of which are power series in $t$ with integer exponents that are bounded below:
\begin{equation}
a_mt^m+a_{m+1}t^{m+1}+a_{m+2}t^{m+2}+\cdots\end{equation}
where $m\in\mathbb{Z}$, $a_i\in K$ for all $i$, and $a_m\neq 0$.
We define a valuation on $K((t))$ by reading off the exponent of the smallest nonzero term:
\begin{equation}
\textrm{val}\left(a_mt^m+a_{m+1}t^{m+1}+a_{m+2}t^{m+2}+\cdots\right)=m.\end{equation}
Show that this is indeed a valuation on $K((t))$.
\end{exercise}

\chproblem{It turns out that the field $K((t))$ is not algebraically closed, even if $K$ is.  For an example of an algebraically closed field with a nontrivial valuation, we turn to the {\emph{field of Puiseux series over $K$}}\index{Puiseux series}, written $K\{\!\{t\}\!\}$.  A nonzero element of this field is of the form

\begin{equation}
a_mt^{m/n}+a_{m+1}t^{{(m+1)}/n}+a_{m+2}t^{{(m+2)}/n}+\cdots\end{equation}
where $m\in\mathbb{Z}$, $n\in\mathbb{Z}^+$, $a_i\in k$ for all $i$, and $a_k\neq 0$.  Note that the value of $n$ can vary between different elements of $K\{\!\{t\}\!\}$, so we could equivalently define a single Puiseux series as a power series in $t$ with rational exponents, where there is a lower bound on the denominator of the exponents.  Again, we can define a valuation by reading off the lowest exponent:
\begin{equation}
\textrm{val}\left(a_mt^{m/n}+a_{m+1}t^{{(m+1)}/n}+a_{m+2}t^{{(m+2)}/n}+\cdots\right)=m/n.\end{equation}
Show that if $K$ is algebraically closed and $\textrm{char}(K)=0$, then $K\{\!\{t\}\!\}$ is algebraically closed.

}

Valuations have a similar flavor to tropical arithmetic, at least if we use the min convention instead of the max convention:  they introduce an infinity element $\infty$, they turn multiplication into addition, and they turn addition into  a minimum (except possibly when the valuations tie).  They also justify the notation of ``vanishing'' as being connected to a minimum or maximum being achieved at least twice, as you'll show in the following exercise.

\begin{exercise} Let $\textrm{val}$ be a valuation on a field $k$, and let $a_1,\ldots,a_n\in k$ with $n\geq 2$.  Show that if $a_1+a_2+\cdots+a_n=0$, then the minimum value among $\textrm{val}(a_1),\ldots,\textrm{val}(a_n)$ occurs at least twice.
\end{exercise}

\subsection{Two Ways To Tropicalize}

To stay consistent with the rest of this chapter, we will continue working in the max convention\footnote{Because we are working in the max convention, there are many instances when we have to consider $-1$ times a valuation.  In the min convention, we can just consider valuations.}.
We now explore two ways of taking a variety $\textrm{V}(I)\subset (k^*)^n$ and moving it into $\mathbb{R}^n$.    One way is be to take coordinate-wise valuation of points in $\textbf{V}(I)$, and append a minus sign onto each coordinate.  That is, we consider the set image of $\textbf{V}(I)$ under the map
\begin{equation}-\textrm{val}:  (k^*)^n\rightarrow \mathbb{R}^n,\end{equation}
\begin{equation}-\textrm{val}(a_1,\cdots,a_n):= (-\textrm{val}(a_1),\cdots,-\textrm{val}(a_1)).\end{equation}
The other way is to consider polynomials $f\in I$, and to turn them into tropical polynomials.  Given $f\in I$ with $f=\sum_\alpha c_\alpha x_1^{\alpha_1}\cdots x_n^{\alpha^n}$, consider the tropical polynomial
\begin{equation}\textrm{trop}(f):=\bigoplus_{\alpha}(-\textrm{val}(c_\alpha))\odot x_1^{\alpha_1}\odot\cdots\odot x_n^{\alpha_n}.\end{equation}
Since $\textbf{V}(I)=\bigcap_{f\in I} \textbf{V}(f)$, we  consider $\bigcap_{f\in I} \mathcal{T}(\textrm{trop}(f))$ as a tropical version of $\textbf{V}(I)$.  We call this intersection the {\emph{tropicalization}}\index{tropicalization} of $\textbf{V}(I)$.

\begin{exercise}
 Let $k=\mathbb{C}\{\!\{t\}\!\}$, and define $f\in k[x,y]$ by
\begin{equation}\begin{matrix}f(x,y)& = &\left(\frac{\sqrt{-1}}{\pi}t^3-3t^{10/3}+\cdots\right)x^2+1000xy\\&&+(1-t^{1/2}+t^{5/8}+\cdots)x+y+(\sqrt{5}t-t^{100}).
\end{matrix}
\end{equation}
Find the tropicalization of $\textbf{V}(f)$.
\end{exercise}

\begin{example}  Let $k=\mathbb{C}\{\!\{t\}\!\}$ where $\mathbb{C}$ is the field of complex numbers, and consider the set $\textbf{V}(I)\subset (k^*)^2$ where $I$ is generated by the single polynomial $x+ty+2\in k[x,y]$.  A point $(a,b)\in \textbf{V}(I)$ is sent to $(-\textrm{val}(a),-\textrm{val}(b))$ by the map $-\textrm{val}$.  Note that if $(a,b)\in \textbf{V}(I)$, then $a=-tb-2$.  This means that either $\textrm{val}(a)=\min\{\textrm{val}(-tb),-2\}=\min\{\textrm{val}(b)+1,0\}$, or $\textrm{val}(a)\geq \min\{\textrm{val}(b)+1,0\}$ with $\textrm{val}(b)+1=0$.  Equivalently, either $-\textrm{val}(a)=\max\{-\textrm{val}(b)-1,0\}$ or $-\textrm{val}(a)\geq \max\{-\textrm{val}(b)-1,0\}$ with $-\textrm{val}(b)-1=0$.  So, all points $(A,B)$ in $-\textrm{val}(\textbf{V}(I))$ fall into one of three classes:
\begin{itemize}
\item  $A=B-1\leq 0$
\item  $A=0\leq B-1$
\item  $B-1=0\leq A$
\end{itemize}  So, the minimum between $A$, $B-1$, and $0$ is achieved at least twice.  In other words, $-\textrm{val}(\textbf{V}(I))\subset \mathcal{T}(x\oplus(-1\odot y)\oplus 0)$.  We do not have equality, since all points in  $-\textrm{val}(\textbf{V}(I))$ have rational coordinates; we leave it as an exercise to show that $-\textrm{val}(\textbf{V}(I))= \mathcal{T}(x\oplus(-1\odot y)\oplus 0)\cap \mathbb{Q}^2$

Note that $\textrm{trop}(x+ty+2)=x\oplus (-1\odot y)\oplus 0$.  All polynomials in $I$ are multiples of  $x+ty+2$, which means that $\bigcap_{f\in I} \mathcal{T}(\textrm{trop}(f))=\mathcal{T}\left(x\oplus (-1\odot y)\oplus 0\right)$.  So, the tropicalization of $\textbf{V}(I)$ is the tropical line defined by $x\oplus (-1\odot y)\oplus 0$.
\end{example}

These two constructions gave us similar, but not identical, subsets of $\mathbb{R}^2$: we had containment of $-\textrm{val}(\textbf{V}(I))$ in the tropicalization of $\textbf{V}(I)$, though these sets were not equal.

\begin{exercise}  Show that we always have $-\textrm{val}(\textbf{V}(I))\subset\bigcap_{f\in I} \mathcal{T}(\textrm{trop}(f))$.
\end{exercise}

It turns out that, as long as we are working over an algebraically closed field, these two sets are equal up to taking a closure in the usual Euclidean topology of $\mathbb{R}^n$.

\begin{theorem}[{The Fundamental Theorem of Tropical Geometry}\index{The Fundamental Theorem of Tropical Geometry}]
Let $k$ be an algebraically closed field with a nontrivial valuation $\textrm{val}$, and let $I$ be an ideal of $k[x_1^{\pm},\cdots,x_n^{\pm1}]$.  Then
\begin{equation}
\overline{-\textrm{val}(\textbf{V}(I))}=\bigcap_{f\in I} \mathcal{T}(\textrm{trop}(f)).\end{equation}
\end{theorem}

This fact is a key result of tropical geometry, originally proved by Kapranov in an unpublished manuscript when $I$ is generated by a single polynomial.  A proof of the more general result appears in \cite[Theorem 3.2.3]{maclagan-sturmfels}.

Given $X=\textbf{V}(I)\subset (k^*)^n$, let $\textrm{Trop}(X)$ denote the set $\overline{-\textrm{val}(\textbf{V}(I))}$.  Understanding the relationship between $X$ and $\textrm{Trop}(X)$ is one of the core themes in tropical geometry.

\subsection{Tropical Intersections}

Let $X$ and $Y$ be varieties in $(k^*)^n$.  Let us consider how $\textrm{Trop}(X\cap Y)$ and $\textrm{Trop}(X)\cap\textrm{Trop}(Y)$ relate to one another.

\begin{exercise}  Show that we always have $\textrm{Trop}(X\cap Y)\subset \textrm{Trop}(X)\cap\textrm{Trop}(Y)$.  (This is mostly an exercise in set theory.)
\end{exercise}

The question then becomes whether we have an equality of these sets.  If we do, then every tropical intersection point in $\textrm{Trop}(X)\cap\textrm{Trop}(Y)$  ``lifts'' to an intersection point in $X\cap Y$.  One core result from \cite{op} is that if $\textrm{Trop}(X)$ and $\textrm{Trop}(Y)$ intersect in components of the expected dimensions, then indeed the points do lift; if $n=2$ and $\textrm{Trop}(X)$ and $\textrm{Trop}(Y)$ are tropical plane curves, this means they intersect in isolated points.  Not only that, these points lift with the expected multiplicity!  If $\textrm{Trop}(X)$ and $\textrm{Trop}(Y)$ intersect in higher dimensional components, the story is more complicated.

\begin{example}\label{example:intersection}  Let $k=\mathbb{C}\{\!\{t\}\!\}$, and let $f,g\in k[x,y]$ be defined by $f(x,y)=ax+by+c$ and $g(x,y)=dx+ey+f$, where $\textrm{val}(a)=\textrm{val}(b)=\textrm{val}(c)=\textrm{val}(d)=\textrm{val}(e)=\textrm{val}(a)=0$.  Let $X=\textbf{V}(f)$ and $Y=\textbf{V}(g)$ be the two lines defined by these equations.  Then $\textrm{Trop}(X)=\textrm{Trop}(Y)=\mathcal{T}(x\oplus y\oplus 0)$, the tropical line in Figure \ref{fig:line}.  This means $\textrm{Trop}(X)\cap \textrm{Trop}(Y)=\mathcal{T}(x\oplus y\oplus 0)$.  Unless $X$ and $Y$ are the same line, at most one of these infinitely many tropical intersection points can lift to an intersection point of $X$ and $Y$.  Let's determine which point might lift.

Assume that $X\cap Y$ consists of one point.  We can solve the equations  $ax+by+c=dx+ey+f=0$ to find the intersection point as $\left(\frac{ce-bf}{bd-ae},\frac{af-cd}{bd-ae}\right)$.  So we know that 
\begin{equation}\textrm{Trop}(X\cap Y)=\textrm{Trop}(\{\left(\frac{ce-bf}{bd-ae},\frac{af-cd}{bd-ae}\right)\})=\{(-\textrm{val}(\frac{ce-bf}{bd-ae}),-\textrm{val}(\frac{af-cd}{bd-ae})\}.\end{equation} If there is no cancellation in $ce-bf,bd-ae,af-cd$, and $bd-ae$, then $\textrm{Trop}(X\cap Y)$ is  $\{(0,0)\}$, which is the stable tropical intersection $\textrm{Trop}(X)\cap_{st}\textrm{Trop}(Y)$.  However, there are cases that give different values for $\textrm{Trop}(X\cap Y)$.  Let $r$ be a positive rational number, and note that:
\begin{itemize}
\item  If $f=x+2y+(1+t^r)$ and $g=x+y+1$, then the intersection point $X\cap Y$ is $\left(-1+t^r,-t^r\right)$, which is sent to $(0,-r)$.
\item  If $f=2x+y+(1+t^r)$ and $g=x+y+1$, then the intersection point $X\cap Y$ is $\left(-t^r,-1+t^r,\right)$, which is sent to $(-r,0)$.
\item  If  $f=(2+t^r)x+2y+1$ and $g=x+y+1$, then the intersection point $X\cap Y$ is $\left(\frac{1}{t^r},\frac{1+t^r}{t^r}\right)=\left(t^{-r},t^{-r}\left(1+t^r\right)\right)$, which is sent to $(r,r)$.
\end{itemize}
This means if all we know about $X$ and $Y$ is that   $\textrm{Trop}(X)\cap \textrm{Trop}(Y)=\mathcal{T}(x\oplus y\oplus 0)$, then \emph{any} point in $\mathcal{T}(x\oplus y\oplus 0)\cap \mathbb{Q}^2$ could be the image of the intersection point of $X$ and $Y$.
\end{example}

\chproblem{Let $a,b,c,d,e,f\in k=\mathbb{C}\{\!\{t\}\!\}$, where $\text{val}(a)=\text{val}(b)=\text{val}(c)=\text{val}(d)=\text{val}(e)=0$ and $\textrm{val}(f)=1$.
Consider the two polynomials $f,g \in k[x,y]$ defined by
\begin{equation}
f(x,y)=ax+by+c,
\end{equation}
\begin{equation}
g(x,y)=dxy+ex+fy.
\end{equation}
Let $X=\textbf{V}(f)$, and $Y=\textbf{V}(g)$.
What are the possible configurations of $\text{Trop}(X\cap Y)$ inside of $\text{Trop}(X)\cap\text{Trop}(Y)$?}

\resproject{Study the possibilities for $\textrm{Trop}(X\cap Y)$ inside of $\textrm{Trop}(X)\cap\textrm{Trop}(Y)$, for plane curves or in higher dimensions.  Some resources to check are \cite{lifting-self, tropical-images,op,or}.}

%
%

\begin{thebibliography}{99.}%
%
%

\bibitem{amm}  Amdeberhan, T., Medina, L. A., Moll, V. H.: Asymptotic valuations of sequences satisfying first order recurrences. Proc. Amer. Math. Soc. \textbf{137},  no. 3 (2009)


\bibitem{blmpr}  Baker, M., Len, Y., Morrison, R.,  Pflueger, N.,  Ren, Q.: Bitangents of tropical plane quartic curves. Math. Z. \textbf{282}, no. 3-4 (2016)

\bibitem{baker-norine} Baker, M., Norine, S.: Riemann-Roch and Abel-Jacobi theory on a finite graph. Adv. Math. 215 no. 2 (2007)

\bibitem{balaban} Balaban, A.T.: Enumeration of cyclic graphs.  Chemical Applications of Graph Theory (A.T. Balaban, ed.) 63?105, Academic Press (1976)

\bibitem{bkt} Bashelor, A., Ksir, A., Traves, W.: Enumerative algebraic geometry of conics. Amer. Math. Monthly \textbf{115}, no. 8  (2008)

\bibitem{bmv} Brannetti, S., Melo, M., Viviani, F.: On the tropical Torelli map. Adv. Math.,
226(3) (2011)

\bibitem{bjms} Brodsky, S., Joswig, M., Morrison, R. Sturmfels, B.: Moduli of tropical plane curves. Res. Math. Sci. \textbf{2}, Art. 4 (2015)

\bibitem{butkovic} Butkovi\u c, P.:  Max-linear systems: theory and algorithms. Springer Monographs in Mathematics. Springer-Verlag London, Ltd., London (2010). 

\bibitem{sprawling}Cartwright, D., Dudzik, A., Manjunath, M., Yao, Y.: Embeddings and immersions of tropical curves. Collect. Math. \textbf{67}, no. 1 (2016)

\bibitem{movingout} Castryck, W.: Moving out the edges of a lattice polygon. Discrete and Computational Geometry \textbf{47}, no. 3 (2012)

\bibitem{cv} Castryck, W., Voight, J.: On nondegeneracy of curves. Algebra and Number Theory \textbf{3} (2009)

\bibitem{chan} Chan, M.: Combinatorics of the tropical Torelli map. Algebra Number Theory,
6(6) (2012)

\bibitem{chan2}  Chan, M.: Tropical hyperelliptic curves. J. Algebraic Combin. \textbf{37}, no. 2 (2013)

\bibitem{k4}  Chan, M., Jiradilok, P.: Theta characteristics of tropical $K_4$-curves. Combinatorial algebraic geometry, 65--86, Fields Inst. Commun., 80, Fields Inst. Res. Math. Sci., Toronto, ON (2017)

\bibitem{cd}  Cools, F., Draisma, J.: On metric graphs with prescribed gonality. J. Combin. Theory Ser. A \textbf{156} (2018)

\bibitem{chip-firing} Corry, S., Perkinson, D: Divisors and sandpiles. An introduction to chip-firing. American Mathematical Society, Providence, RI (2018)



\bibitem{clo} Cox, D. A., Little, J., O'Shea, D.: Ideals, varieties, and algorithms. An introduction to computational algebraic geometry and commutative algebra. Fourth edition. Undergraduate Texts in Mathematics. Springer, Cham, xvi+646 pp. (2015)

\bibitem{de} de Bruijn, N. G., Erdös, P: On a combinatorial problem. Nederl. Akad. Wetensch., Proc. \textbf{51} (1948) 

\bibitem{triangulations}  De Loera, J.A., Rambau, J., Santos, F.: Triangulations. Structures for algorithms and applications. Algorithms and Computation in Mathematics, 25. Springer-Verlag, Berlin (2010)

\bibitem{sevencircles}  Evelyn, C. J. A., Money-Coutts, G. B.,Tyrrell, J. A.: The seven circles theorem and other new theorems. Stacey International, London (1974)

\bibitem{fno}  Farouki, R.T., Neff, C., O'Conner, M.A.:  Automatic parsing of degenerate quadric-surface intersections.  ACM Transactions on Graphics (TOG) \textbf{8}, No. 3 (1989)

\bibitem{gk}  Gathmann, A., Kerber, M.: A Riemann-Roch theorem in tropical geometry. Math. Z. \textbf{259}, no. 1 (2008)

\bibitem{polymake}   Gawrilow, E., Joswig, M.: polymake: a framework for analyzing convex polytopes. Polytopes-combinatorics and computation (Oberwolfach, 1997), 43?73, DMV Sem., 29, Birkhäuser, Basel (2000)

\bibitem{tropicalschemes}  Giansiracusa, J., Giansiracusa, N.: Equations of tropical varieties. Duke Math. J. \textbf{165}, no. 18, 3379--3433 (2016)

\bibitem{m2} Grayson, D., Stillman, M. E.:  Macaulay2, a software system for research in algebraic geometry. Available at \url{http://www.math.uiuc.edu/Macaulay2/}



\bibitem{hmrt}  Hahn, M.A.,  Markwig, H.,  Ren, Y., Tyomkin, I.: Tropicalized quartics and canonical embeddings for tropical curves of genus 3. arXiv preprint \texttt{arXiv:1802.02440} (2018)



\bibitem{hartshorne}
Hartshorne, R.: Algebraic geometry. Graduate Texts in Mathematics, No. 52. Springer-Verlag, New York-Heidelberg {\rm xvi}+496 pp. (1977)

\bibitem{hilbert}
Hilbert, D.: Ueber die Theorie der algebraischen Formen. (German) Math. Ann. \textbf{36},  no. 4 (1890)

\bibitem{gfan}   Jensen, A.: Gfan, a software system for Gr\"{o}bner fans and tropical varieties. Available at
\url{http://home.imf.au.dk/jensen/software/gfan/gfan.html}.

\bibitem{kz} Kaibel, V., Ziegler, G.M.: Counting lattice triangulations. Surveys in combinatorics, 2003 (Bangor), 277--307, London Math. Soc. Lecture Note Ser., 307, Cambridge Univ. Press, Cambridge, (2003)

\bibitem{koelman}   Koelman, R.: The number of moduli of families of curves on toric surfaces.
Ph.D. thesis, Katholieke Universiteit Nijmegen (1991)

\bibitem{lz}  Lagarias, J.C., Ziegler,  G.M.: Bounds for lattice polytopes containing a fixed number of
interior points in a sublattice. Canadian J. Math. \textbf{43} (1991)

\bibitem{legall} Le Gall, François: Powers of tensors and fast matrix multiplication. Proceedings of the 39th International Symposium on Symbolic and Algebraic Computation (ISSAC 2014)



\bibitem{lifting}  Len, Y., Markwig, H..: Lifting tropical bitangents. arXiv preprint \texttt{arXiv:1708.04480} (2018)

\bibitem{lifting-self}  Len, Y., Satriano., M.: Lifting tropical self intersections. arXiv preprint \texttt{arXiv:1806.01334} (2018)


\bibitem{lengyel} Lengyel, T.: On the divisibility by $2$ of the Stirling numbers of the second kind. Fibonacci Quart. \textbf{32}, no. 3 (1994)


\bibitem{lin-tran}  Lin, B., Tran, N. M.: Linear and rational factorization of tropical polynomials. arXiv preprint \texttt{arXiv:1707.03332} (2017)

\bibitem{tropicalideals} Maclagan, D., Rinc\'{o}n, F.: Tropical ideals. Compos. Math. \textbf{154}, no. 3   (2018)

\bibitem{maclagan-sturmfels}  Maclagan, D., Sturmfels, B.: Introduction to tropical geometry. Graduate Studies in Mathematics, 161. American Mathematical Society, Providence, RI (2015)

\bibitem{fibonacci}  Medina, L. A., Rowland, E.: $p$-regularity of the $p$-adic valuation of the Fibonacci sequence. Fibonacci Quart. \textbf{53}, no. 3 (2015)

\bibitem{mz}  Mikhalkin, G., Zharkov, I.: Tropical curves, their Jacobians and theta functions. Curves and abelian varieties, 203--230, Contemp. Math., 465, Amer. Math. Soc., Providence, RI (2008)


\bibitem{morrison-hyperelliptic}  Morrison, R.: Tropical hyperelliptic curves in the plane. arXiv preprint \texttt{arXiv:1708.00571} (2017)

\bibitem{tropical-images}  Morrison, R: Tropical images of intersection points. Collect. Math. \textbf{66}, no. 2 (2015)


\bibitem{op} Osserman, B., Payne, S.: Lifting tropical intersections. Doc. Math. \textbf{18} (2013)

\bibitem{or}  Osserman, B., Rabinoff, J.: Lifting nonproper tropical intersections. Tropical and non-Archimedean geometry, 15--44, Contemp. Math., 605, Centre Rech. Math. Proc., Amer. Math. Soc., Providence, RI, (2013)


\bibitem{pick} Pick, G. A.:  Geometrisches zur Zahlenlehre, Sitzenber. Lotos (Prague) \textbf{19} (1899)

\bibitem{pin} Pin, J.: Tropical semirings. Idempotency (Bristol, 1994), 50--69, Publ. Newton Inst., 11, Cambridge Univ. Press, Cambridge (1998)


\bibitem{plucker}  Pl\"{u}cker, J.: Solution d'une question fondamentale concernant la théorie générale des courbes. (French) J. Reine Angew. Math. \textbf{12} (1834)


\bibitem{firststeps}  Richter-Gebert, J., Sturmfels, B., Theobald, T.: First steps in tropical geometry. Idempotent mathematics and mathematical physics, 289--317, Contemp. Math., \textbf{377}, Amer. Math. Soc., Providence, RI (2005)

\bibitem{topcom} Rambau, J.: TOPCOM: Triangulations of Point Configurations and Oriented Matroids, Mathematical Software - ICMS 2002 (Cohen, Arjeh M. and Gao, Xiao-Shan and Takayama, Nobuki, eds.), World Scientific, pp. 330-340 (2002)

\bibitem{scott}  Scott, P.R.: On convex lattice polygons. Bull. Austral. Math. Soc. \textbf{15} (1976)

\bibitem{simon}  Simon, I.: Recognizable sets with multiplicities in the tropical semiring. Mathematical foundations of computer science, 1988 (Carlsbad, 1988), 107--120, Lecture Notes in Comput. Sci., 324, Springer, Berlin, (1988)

\bibitem{strassen}  Strassen, V.: Gaussian Elimination is not Optimal. Numer. Math. \textbf{13} (1969)

\end{thebibliography}
%

\printindex
\end{document}